\documentclass[12pt]{article}
\usepackage{amsmath,amsfonts,amsthm,amssymb,amscd}
\renewcommand{\deg}{\mathop{\rm deg}\nolimits}

\newcommand{\C}{{\mathbb C}}
\newcommand{\R}{{\mathbb R}}

\newcommand{\Z}{{\mathbb Z}}

\newcommand{\N}{{\mathbb N}}

\newcommand{\supp}{\mathop{\rm supp}\nolimits}

\theoremstyle{plain}
\newtheorem{theorem}{Theorem}[section]
\newtheorem{lemma}[theorem]{Lemma}
\newtheorem{proposition}[theorem]{Proposition}

\theoremstyle{definition}

\theoremstyle{remark}
\newtheorem{remark}[theorem]{Remark}

\def\a{\alpha}

\numberwithin{equation}{section}

\def\supp{\operatorname{supp}}

\def\im{\operatorname{im}}
\def\re{\operatorname{re}}

\def\a{\alpha}

\def\l{\lambda}

\def\be{\begin{equation}}
\def\ee{\end{equation}}
\begin{document}
\author{Galina Perelman\footnote{LAMA, UMR CNRS 8050, Universit\'e Paris-Est Cr\'eteil, 
61, avenue du G\'en\'eral de Gaulle, 94010
Cr\'eteil Cedex, France}}
\title{Blow up dynamics for equivariant critical Schr\"odinger maps}
%\date{10.08}
\date{}
%\date{(preliminary version)}
\maketitle

\begin{abstract}

For the Schr\"odinger map equation 
$u_t=u\times \triangle u$
in
$\R^{2+1}$, with values in $S^2$,  we prove for any $\nu>1$ 
the existence of equivariant finite time blow up solutions
of the form $u(x,t)=\phi(\lambda(t) x)+\zeta(x,t)$,
where $\phi$ is a lowest energy steady state, $\lambda(t)=t^{-1/2-\nu}$
and $\zeta(t)$ is arbitrary small in $\dot H^1 \cap \dot H^2$.

%the parameter $\nu>1/2$ can be chosen arbitrarily.

\end{abstract}

\setcounter{section}{0}

\section{Introduction}
\label{s1}

\subsection{Setting of the problem and statement of the result}

In this paper we consider the Schr\"odinger flow for maps from $\R^2$ to $S^2$:
\begin{equation}\label{0.1}\begin{split}
&u_t=u\times \triangle u, \quad x=(x_1, x_2)\in \R^2,\,\, t\in\R,\\
&u|_{t=0}=u_0,\end{split}
\end{equation}
where $u(x,t)=(u_1(x,t), u_2(x,t), u_3(x,t))\in S^2\subset \R^3$.

Equation \eqref{0.1} conserves the energy
\begin{equation}\label{0.2}
E(u)={1\over 2}\int_{\R^2}dx |\nabla u|^2.
\ee
The problem is critical in the sense that both \eqref{0.1} and \eqref{0.2}
are invariant with respect to the scaling
$u(x,t)\rightarrow u(\lambda x, \lambda^2 t)$, $\lambda\in \R_+$.
\par To a finite energy map $u: \R^2\rightarrow S^2$ one can associate the degree:
$$\deg(u)={1\over 4\pi}\int_{\R^2}dx u_{x_1}\cdot J_uu_{x_2}, $$
where $J_u$ is defined by 
$$J_uv=u\times v,\quad v\in \R^3.$$
It follows  from \eqref{0.2} that
\begin{equation}\label{0.3}
E(u)\geq 4\pi|\deg(u)|.
\ee
This inequality is saturated by the harmonic maps $\phi_m$, $m\in \Z^+$:
\be\label{0.4}\begin{split}
&\phi_m(x)=e^{m\theta R}Q^m(r), \quad Q^m=(h_1^m, 0, h_3^m)\in S^2,\\\
&h_1^m(r)=\frac{2r^m}{r^{2m}+1},\quad h_3^m(r)=\frac{r^{2m}-1}{r^{2m}+1},\end{split}
\ee
Here $(r,\theta)$ are polar coordinates in  $\R^2$: $x_1+ix_2=e^{i\theta}r$, 
and $R$ is the generator of the horizontal rotations:
$$ R=\left(\begin{array}{ccc}
0&-1&0\\
1&0&0\\
0&0&0\end{array}\right) ,$$
or equivalently
$$Ru= k\times u, \quad k=(0,0,1).$$
One has 
$$\deg \phi_m=m, \quad E(\phi_m)=4\pi m.$$
Up to the symmetries $\phi_m$ are  the only energy minimizers in their homotopy class.

Since $\phi_1$ will play a central role in the analysis developed in this paper, we set
$\phi=\phi_1$, $Q=Q_1$, $h_1=h_1^1$, $h_3=h_3^1$.

The local/global well-posedness of \eqref{0.1} has been extensively studied in past years.
Local existence for smooth initial data goes back to \cite{SSB}, see also 
\cite{McG}.
The case of small data of low regularity was studied in several works, 
the definite result being obtained by Bejenaru, Kenig and Tataru in \cite{BKT},
where the global existence and scattering was proved for general
$\dot H^1$ small initial data.
Global existence for equivariant small energy initial data was proved earlier in \cite{CHU} (by $m$-equivariant map $u:\R^2\rightarrow S^2\subset \R^3$, $m\in\Z^+$
one means  a map of the form
$u(x)=e^{m\theta R}v(r),$
where $v:\R_+\rightarrow S^2\subset\R^3$,
$m$-equivariance being preserved by the Schr\"odinger flow \eqref{0.1}).
In the radial case $m=0$, the global existence for $H^2$ data was established by Gustafson and Koo \cite{GK}.
Very recently,  Bejenaru, Ionescu, Kenig and Tataru \cite{BIKT} proved global existence and scattering for equivariant  data with energy less than $4\pi$.
The dynamics of $m$-equivariant Schr\"odinger maps 
with initial data close to
$\phi_m$ was studied by Gustafson, Kang,  Nakanishi, and
Tsai \cite{GKT1}, \cite{GKT2},  \cite{GNT}
and later  by Bejenaru,  Tataru \cite{BT} in the case $m=1$. 
The stability/instability results of these works strongly suggest a possibility
of regularity breakdown in solutions of \eqref{0.1}
via concentration of the lowest energy harmonic map $\phi$.
For a closely related model of wave maps
this type of regularity breakdown was  proved 
by Kriger, Schlag and Tataru \cite{KST1} and by  Raphael and Rodnianski \cite{RR}.
These authors showed the existence of 1- equivariant blow up solutions close to 
$\phi(\lambda(t)x)$ with $\l(t)\sim \frac{e^{\sqrt{|\ln(T^*-t)|}}}{T^*-t}$ as $t\rightarrow T^*$
\cite{RR}, and with $\l(t)\sim \frac{1}{(T^*-t)^{1+\nu}}$ as $t\rightarrow T^*$ where $\nu>1/2$
can be chosen arbitrarily \cite{KST1}
(here $T^*$ is the blow up time).  While the blow up dymamics exhibited in \cite {RR} is stable (in some
strong topology), the continuum of blow up solutions constructed by Kriger, Schlag and Tataru is believed to be non-generic. Recently, the results of \cite{RR} were generalized to the case of Schr\"odinger map equation \eqref{0.1} by Merle, Raphael and Rodnianski in \cite{MRR} where they proved
 the existence of 1-equivariant blow up solutions of \eqref{0.1}
close to $\phi(\lambda(t)x)$ with $\l(t)\sim \frac{(\ln(T^*-t))^2}{T^*-t}$.

Our objective in this paper is to show that \eqref{0.1} also admits 1-equivariant
Kriger-Schlag-Tataru type blow up solutions
that correspond
to certain initial data of the form
$$u_0=\phi+\zeta_0,$$
where $\zeta_0$ is 1-equivariant and can be chosen arbitrarily small in $\dot H^1 \cap \dot H^3$.
Let us recall   (see \cite{GKT1},  \cite{GKT2}, \cite{GNT}, \cite{BT}) that such initial data result in unique local solution of the same
regularity, and  as long as the solution exists it
 stays $\dot H^1$ close to a two parameter family of  1-equivariant harmonic maps 
$\phi^{\a,\l}$ , $\a\in \R/2\pi \Z$, $\l\in\R_+$ generated from $\phi$
by rotations and scaling:
$$\phi^{\a,\l}(r,\theta)=e^{\a R}\phi(\l r,\theta).$$
%Furthermore,
%the length scale $\l(t)$ of the $\dot H^1$ nearest harmonic map satisfies
%$$\l(t)\geq C\|u(t)\|_{\dot H^2}.$$
%see ... below for a more precise statement.

The following theorem is the main result of this paper.

\begin{theorem}
For any $\nu>{1}$, $\alpha_0\in \R$, and any $\delta >0$ sufficiently small
there exist $t_0>0$ and a 1-equivariant solution $ u\in C((0,t_0], \dot H^1\cap \dot H^{3})$ of \eqref{0.1}
of the form:
\be\label{0.10}
u(x,t)=e^{\alpha(t)R}\phi(\lambda(t)x)+\zeta(x,t),\ee
where 
\be\label{0.11}
\lambda(t)=t^{-1/2-\nu},\quad \alpha(t)=\alpha_0\ln t,\ee
\be\label{0.12}
\|\zeta(t)\|_{\dot H^1\cap \dot H^{2}}\leq \delta, 
\quad \|\zeta(t)\|_{\dot H^3}\leq C_{\nu,\alpha_0}t^{-1},
\quad\forall t\in (0,t_0].\ee
Furthermore, as $t\rightarrow 0$, $\zeta(t)\rightarrow \zeta^*$ in $\dot H^1\cap \dot H^{2}$ 
with $\zeta^*\in H^{1+2\nu -}$.
\end{theorem}

\begin{remark} In fact, using the arguments developed in this paper one can show that
the same result remains valid with $\dot H^3$ replaced by $\dot H^{1+2s}$ for any 
$ 1\leq s<\nu$.
\end{remark}

\subsection{Strategy of the proof} The proof of Theorem 1.1 contains two main steps. The first step is a construction of  an approximate solutions $u^{(N)}$ that have form \eqref{0.10},\eqref{0.11},
\eqref{0.12},
 and solve \eqref{0.1} up an arbitrarily high order error $O(t^N)$, very much in the spirit
of the work of of Kriger, Schlag and Tataru
 \cite{KST1}.
\par The second step is to  build the exact solution
by solving the problem for the small remainder forward in time with zero initial data at $t=0$.
The control of remainder is  achieved by means of suitable energy type estimates,
see section 3 for the detail.
\section{Approximate solutions}
\subsection{Preliminaries}
We consider \eqref{0.1} under the 1-equivariance assumption
\be\label{1.001}
u(x,t)=e^{\theta R}v(r,t), \quad v=(v_1,v_2, v_3)\in S^2\subset \R^3.\ee
Restricted to the 1-equivariant functions \eqref{0.1} takes the form
\be\label{1.1}
v_t=v\times (\triangle v+{R^2\over r^2}v),\ee

the energy being given by
$$E(u)=\pi\int_0^\infty dr r(|v_r|^2+{v_1^2+v_2^2\over r^2}).$$
$Q=(h_1, 0, h_3)$ is a stationary solution of \eqref{1.1} and one has the relations
\be\label{1.*2*}
\partial_r h_1=-\frac{h_1h_3}{r},\quad \partial_r h_3=\frac{h_1^2}{r},\ee
\be\label{1.3}
\triangle Q +{R^2\over r^2}Q=\kappa(r)Q,\quad \kappa(r)=-{2h_1^2\over r^2}.\ee

\par The goal of the present section is to prove the following result.
\begin{proposition}\label{p1}
For  any $\delta>0$ sufficiently small and any $N$ sufficiently large  there exists an approximate solution $u^{(N)}:\R^2\times\R_+^*\rightarrow S^2$ of \eqref{0.1} such that the following holds.\\
\noindent (i) $u^{(N)}$ is a $C^\infty$ 1-equivariant profile of the form:
$u^{(N)}=e^{\a(t)R}(\phi(\lambda(t)x)+\chi^{(N)}(\lambda(t)x,t))$,
where $\chi^{(N)}(y,t)= 
e^{\theta R}Z^{(N)}(\rho,t)$,  $\rho=|y|$, verifies
\begin{align}
&\|\partial_\rho Z^{(N)}(t)\|_{L^2(\rho d\rho)}, \|\rho^{-1}Z^{(N)}(t)\|_{L^2(\rho d\rho)}, \|\rho\partial_\rho Z^{(N)}(t)\|_\infty\leq C\delta^{2\nu},\label{p1.1.1}\\
&\|\rho^{-l}\partial_\rho^kZ^{(N)}(t)\|_{L^2(\rho d\rho)}\leq C\delta^{2\nu-1}t^{1/2+\nu},\quad k+l=2,
\label{p1.1.2}\\
&\|\rho^{-l}\partial_\rho^kZ^{(N)}(t)\|_{L^2(\rho d\rho)}\leq Ct^{2\nu},\quad k+l=3,
\label{p*1.1.2}\\
&\|\partial_\rho Z^{(N)}(t)\|_\infty, \|\rho^{-1}Z^{(N)}(t)\|_\infty\leq C\delta^{2\nu-1}t^{\nu},
\label{p1.1.3}\\
&\|\rho^{-l}\partial_\rho^kZ^{(N)}(t)\|_\infty\leq C t^{2\nu},\quad 2\leq l+k\leq 3,\label{p1.1.4}
\end{align}
for any $0<t\leq T(N,\delta)$ with some $T(N,\delta)>0$.  The constants $C$ here and below are independent of $N$ and $\delta$.\\
In addition, one has 
\be\label{p1.1.6}
\|\chi^{(N)}(t)\|_{\dot W^{4,\infty}}+\|<y>^{-1}\chi^{(N)}(t)\|_{\dot W^{5,\infty}}\leq Ct^{2\nu},
\ee
and $<x>^{2(\nu-1)}\nabla^4u^{(N)}(t),\, <x>^{2(\nu-1)}\nabla^2u^{(N)}_t(t)\in L^\infty(\R^2)$.

Furthermore, there exists $\zeta^*_N\in \dot H^{1}\cap  \dot H^{1+2\nu-}$ such that
as $t\rightarrow 0$,\\ $e^{\alpha(t)R}\chi^{(N)}(\lambda(t)\cdot, t)\rightarrow \zeta^*_N$ in
$\dot H^{1}\cap  \dot H^{2}$.\\
%Furthermore, $v^{(N)}(t)-1\in <x>^{-1}L^2(rdr)$ and
\noindent (ii) The corresponding error $r^{(N)}=-u^{(N)}_t+u^{(N)}\times \Delta u^{(N)}$
verifies
\be\label{p1.2}
\|{ r}^{(N)}(t)\|_{H^{3}}+ \|\partial_t{ r}^{(N)}(t)\|_{H^{1}}+
\|<x>r^{(N)}(t)\|_{L^{2}}
\leq t^N,
\quad 0<t\leq T(\delta,N).\ee
\end{proposition}

\noindent{\it Remarks}.\\
1. Note that estimates \eqref{p1.1.1}, \eqref{p1.1.2} imply:
\be\label{p1.11}
\|u^{(N)}(t)-e^{\a(t)R}\phi(\l(t)\cdot)\|_{\dot H^1\cap \dot H^{2}}\leq \delta^{2\nu -1} , \quad\forall t\in (0,T(N,\delta)].\ee
2. It follows from our construction that $\chi^{(N)}(t)\in \dot H^{1+2s}$ for any $ s<\nu$
with the estimate
$ \|\chi^{(N)}(t)\|_{\dot H^{1+2s}(\R^2)}\leq C(t^{2\nu}+ t^{s(1+2\nu)}\delta^{2\nu-2s})$.\\
3. The remainder $r^{(N)}$ verifies in fact, for any $m, l ,k$,
$$\|<x>^l\partial_t^m{ r}^{(N)}(t)\|_{H^{k}}\leq C_{l,m,k}t^{N-C_{l,m,k}},$$
provided $N\geq C_{l,m,k}$.

We will give  the proof of  proposition \ref{p1} in the case of $\nu $ irrational only, which allows to slightly simplify the presentation.
The extension to $\nu $ rational  is straightforward.

To construct an arbitrarily good approximate solution we analyze separately the three regions that
correspond to three different space scales:
the inner region with the scale $r\lambda(t)\lesssim1$, the self-similar region
where $r=O(t^{1/2})$, and finally the remote region where $r=O(1)$. The inner region is the region where
the blowup concentrates. In this region the solution will be constructed as a perturbation of the profile
$e^{\a(t)R}Q(\l(t)r)$. The self-similar and remote regions are the regions where the solution is  close to $k$
and is  described essentially by the corresponding linearized equation. In the self-similar region the profile of the solution will be determined uniquely by the matching conditions coming out of the inner region, while in the remote region 
the profile remains essentially a free parameter of the construction, only the limiting behavior at the origin is prescribed by the matching process, see subsections 2.3 and 2.4 for the details, see  also
\cite{AH}, \cite{Ber1} for some
 closely related considerations in the context of
the critical harmonic map heat flow.

\subsection{Inner region  $r\lambda(t)\lesssim1$}

We start by  considering  the inner region
$0\leq r\lambda(t)\leq 10t^{-\nu+\varepsilon_1}$,
where $0<\varepsilon_1<\nu$ to be fixed later.
Writing $v(r,t)$ as 
$$v(r,t)=e^{\alpha(t)R}V(\lambda (t)r, t), \quad V=(V_1, V_2, V_3),$$
we get from \eqref{1.1}
\be\label{1.2}
t^{1+2\nu}V_t+\alpha_0t^{2\nu}RV-(\nu+{1\over 2})\rho V_\rho=
V\times(\triangle V+{R^2\over \rho^2}V),\quad \rho=\l(t)r.\ee
We look for a solution of \eqref{1.2} as a perturbation of the harmonic map 
profile $Q(\rho)$. Write
$$V=Q+Z,$$
and further decompose
$Z$ as
$$Z(\rho,t)=z_1(\rho, t)f_1(\rho)+z_2(\rho, t)f_2(\rho)+\gamma(\rho,t)Q(\rho),
$$
where $f_1$, $f_2$ is the orthonormal frame on $T_QS^2$ given by
$$f_1(\rho)=\left(\begin{array}{c}
h_3(\rho)\\
0\\
-h_1(\rho)
\end{array}\right) 
\quad f_2(\rho)=\left(\begin{array}{c}
0\\
1\\
0
\end{array}\right) .$$
On has $$\gamma=\sqrt{1-|z|^2}-1=O(|z|^2),\quad z=z_1+iz_2.$$
Note also the relations 
\begin{equation*}\begin{split}
&\partial_\rho Q=-{h_1\over \rho}f_1,\quad \partial_\rho f_1={h_1\over \rho}Q,
\quad f_2=Q\times f_1,\\
&\Delta f_1+\frac{R^2}{\rho^2}f_1=-\frac{1}{\rho^2}f_1-\frac{2h_3h_1}{\rho^2}Q.
\end{split}
\end{equation*}

Rewriting \eqref{1.2} in terms of $z$ one gets:
\be\label{1.4}\begin{split}
&it^{1+2\nu}z_t-\a_0t^{2\nu}h_3z-i({1\over 2}+\nu)t^{2\nu}\rho z_\rho
=Lz+F(z)+dt^{2\nu}h_1,\\
&L=-\triangle +{1-2h_1^2\over \rho^2},\quad d=\a_0-i(\frac12+\nu),\\
&F(z)=\gamma Lz+z(\triangle\gamma+ {2 h_1\over \rho}\partial_\rho z_1-{2h_1h_3z_1\over \rho^2})
+\frac{2h_1}{\rho}(1+\gamma)\gamma_\rho+dt^{2\nu}\gamma h_1.\\
\end{split}
\ee
Note that  $F$ is at least quadratic in $z$.

We look for a solution of \eqref{1.4} as a power expansion in $t^{2\nu}$:
\be\label{1.5}
z(\rho,t)=\sum_{k\geq 1} t^{2\nu k}z^k(\rho).\ee
Substituting \eqref{1.5} into \eqref{1.4} we get the following
 recurrent system for $z^k,\, k\geq 1$:
\be\label{1.6}
Lz^1=dh_1,
\ee
\be\label{1.7}
Lz^k={\cal F}_k,\quad k\geq 2,
\ee
where ${\cal F}_k$ depends on $z^j,\,\, j=1,\dots, k-1$ only. 
We subject \eqref{1.6}, \eqref{1.7} to zero initial conditions at $\rho=0$:
\be\label{1.8}
z^k(0)=\partial_\rho z^k(0)=0.\ee
\begin{lemma}\label{l1}
 System \eqref{1.6}, \eqref{1.7},  \eqref{1.8}  has a unique solution $(z^k)_{k\geq 1}$, with $z^k \in C^\infty(\R_+)$ for all $k\geq 1$.
In addition, one has:\\
\noindent (i)  $z^k$ has an  odd Taylor expansion  at $0$ that starts at oder $2k+1$;\\
\noindent (ii) as $\rho\rightarrow \infty$, $z^k$ has the following asymptotic expansion
\be\label{1.9}
z^k(\rho)=\sum\limits_{l=0}^{2k}\sum\limits_{j\leq k-(l-1)/2}c^k_{j,l}\rho^{2j-1}(\ln \rho)^l,\ee
with some constants $c^k_{j,l}$.
The asymptotic expansion \eqref{1.9} can be differentiated any number of times with respect to $\rho$.

\end{lemma}

{\it Proof}.
First note that the equation $Lf=0$ has two explicit solutions: $h_1(\rho)$ and
$h_2(\rho)={\rho^4+4\rho^2\ln \rho-1\over \rho(\rho^2+1)}$.

Consider the case $k=1$:
$$Lz^1=dh_1,$$
$$
z^1(0)=\partial_\rho z^1(0)=0.$$
One has
\be\label{1.10}\begin{split}
&z^1(\rho)={d\over 4}\int_0^\rho ds s( h_1(\rho)h_2(s)-h_1(s)h_2(\rho))h_1(s)\\
&={d\rho\over (1+\rho^2)}\int_0^\rho ds {s(s^4+4s^2\ln s -1)\over (1+s^2)^2} - {d(\rho^4+4\rho^2\ln \rho-1)\over \rho(\rho^2+1)}\int_0^\rho ds {s^3\over (1+s^2)^2}\end{split}\ee
Since $h_1$ is a $C^\infty$ function that has an odd Taylor expansion at $\rho=0$ with a linear leading term,
one can easily write an odd Taylor series for $z^1$  with a cubic leading  term, which proves (i) for $k=1$.

The asymptotic behavior of $z^1$ at infinity can be obtained directly from the representation \eqref{1.10}.
As claimed, one has
$$ z^{1}(\rho)=c_{1,0}^1\rho +c_{1,1}^1\rho\ln \rho  +\sum_{j\leq 0}\sum_{l=0,1,2} c_{j,l}^1\rho^{2j-1}(\ln \rho)^l,$$
with $c_{1,0}^1=-c_{1,1}^1=d_1$.

Consider $k>1$. Assume that $z^j, \, j\leq  k-1$, verify (i) and (ii). Then, using \eqref{1.4}, one can easily check that
${\cal F}_k$ is an odd  $C^\infty$  function vanishing at $\rho=0$ at  order $2k-1$, with the following asymptotic expansion as $\rho \rightarrow \infty$:
\begin{equation*}\begin{split}
{\cal F}_k&=\sum_{j=1}^{k-1}\sum_{l=0}^{2k-2j-1}\alpha_{j,l}^k \rho^{2j-1}(\ln \rho)^l+
\sum_{l=0}^{2k-2}\alpha_{0,l}^k \rho^{-1}(\ln \rho)^l\\
&+\sum_{l=0}^{2k-1}\alpha_{-1,l}^k \rho^{-3}(\ln \rho)^l+
\sum_{j\leq -2}\sum_{l=0}^{2k}\alpha_{j,l}^k \rho^{2j-1}(\ln \rho)^l.\end{split}\end{equation*}
As a consequence, $z^k(\rho)={1\over 4}\int_0^\rho ds s( h_1(\rho)h_2(s)-h_1(s)h_2(\rho)){\cal F}_k(s)$
is a $C^\infty$ function with an odd Taylor series at zero starting at order $2k+1$ and as
$\rho\rightarrow \infty$,
$$z^k(\rho)=\sum\limits_{l=0}^{2k}\sum\limits_{j\leq k-(l-1)/2}c^k_{j,l}(\ln \rho)^l\rho^{2j-1},$$
as required. This concludes the proof of lemma \ref{l1}.$\quad\quad\quad\quad\Box$

Returning to $v$ we get a formal solution of \eqref{1.1} of the form
\be\label{1.200}
v(r,t)=e^{\a(t)R}V(\l(t)r,t), \quad V(\rho, t)=Q+\sum_{k\geq 1} t^{2\nu k}Z^k(\rho), \ee
$ Z^k=(Z^k_1, Z^k_2, Z^k_3)$,
where $Z^k_i$, $i=1,2,$  are smooth odd functions of $\rho$ vanishing at $0$ at order $2k+1$,
and $Z^k_3$ is an even function vanishing at  zero at oder $2k+2$. As $\rho \rightarrow \infty$, one has
\be\label{1.20}\begin{split}
&Z_i^k(\rho)=\sum\limits_{l=0}^{2k}\sum\limits_{j\leq k-(l-1)/2}c^{k,i}_{j,l}(\ln \rho)^l\rho^{2j-1},\quad i=1,2,\\
&Z_3^k(\rho)=\sum\limits_{l=0}^{2k}\sum\limits_{j\leq k+1-l/2}c^{k,3}_{j,l}(\ln \rho)^l\rho^{2j-2},\end{split}\ee
with some real coefficients $c^{k,i}_{j,l}$ verifying
$$c^{k,3}_{k+1,0}=0,\quad\forall k\geq 1.$$
The asymptotic expansions \eqref{1.20} can be differentiated any number of times with respect to $\rho$.

Note that in the limit $\rho\rightarrow \infty$, $y\equiv rt^{-1/2}\rightarrow  0$, expansion 
 \eqref{1.200}, \eqref{1.20}
give at least formally
\be\label{1.21}\begin{split}
&V_i(\lambda(t)r,t)=
\sum_{j\geq 0}t^{\nu(2j+1)}\sum_{l=0}^{2j+1}(\ln y -\nu\ln t)^lV_i^{j,l}(y),\quad i=1,2,
\\
&V_3(\lambda(t)r,t)=1+
\sum_{j\geq 1}t^{2\nu j}\sum_{l=0}^{2j}(\ln y -\nu\ln t)^lV_3^{j,l}(y),
\\
&V_i^{j,l}(y)=\sum_{k\geq -j+l/2}c^{k+j,i}_{k,l}y^{2k-1},\quad i=1,2,\\
&V_3^{j,l}(y)=\sum_{k\geq -j+l/2}c^{k+j,3}_{k+1,l}y^{2k},
\end{split}\ee
where the coefficients $c^{k,i}_{j,l}$ with $k\neq 0$ are defined by
\eqref{1.20} and $c^{0,i}_{j,0}$ come from the expansion of $Q$ as $\rho\rightarrow \infty$:
$$h_1(\rho)=\sum\limits_{j\leq 0}c^{0,1}_{j,0}\rho^{2j-1},
\quad h_3(\rho)=1+\sum\limits_{j\leq 0}c^{0,3}_{j,0}\rho^{2j-2},\quad
c^{0,2}_{j,0}=0.$$

For $N\geq2$ define
$$z^{(N)}_{in}=\sum_{k=1}^N z^k,\quad z^{(N)}_{in}=z_{in,1}^{(N)}+iz_{in,2}^{(N)}.$$
Then $z_{in}^{(N)}$ solves \eqref{1.4} up to the error
$X_N=-it^{1+2\nu}\partial_tz^{(N)}_{in}-\a_0t^{2\nu}h_3z^{(N)}_{in}+i({1\over 2}+\nu)t^{2\nu}\rho \partial_\rho z^{(N)}_{in}
+dt^{2\nu}h_1
+Lz^{(N)}_{in}+F(z^{(N)}_{in}),$
verifying
\be\label{error: 1}
|\rho^{-l}\partial^k_\rho \partial_t^mX_N|\leq C_{k,l,m}t^{2\nu N-m}<\rho>^{2N-1-l-k}\ln(2+\rho),\ee
for any $k,m\in\N,\, 0\leq l\leq (2N+1-k)_+$,
$0\leq \rho\leq10 t^{-\nu+\varepsilon_1}$,
$0<t\leq T(N)$, with some  $T(N)>0$.

Set
$$\gamma^{(N)}_{in}=\sqrt{1-|z^{(N)}_{in}|^2}-1,$$
$$Z_{in}^{(N)}=z_{in,1}^{(N)}f_1+z_{in,2}^{(N)}f_2+\gamma^{(N)}_{in}Q,$$
$$V_{in}^{(N)}=Q+Z_{in}^{(N)}\in S^2.$$
Then $V_{in}^{(N)}$ solves
\be\label{eq: V^N}
t^{1+2\nu}\partial_tV_{in}^{(N)}+\alpha_0t^{2\nu}RV_{in}^{(N)}-
(\nu+{1\over 2})\rho \partial_\rho V_{in}^{(N)}=
V_{in}^{(N)}\times(\triangle_\rho V_{in}^{(N)}+{R^2\over \rho^2}V_{in}^{(N)})+{\cal R}_{in}^{(N)},
\ee with
${\cal R}_{in}^{(N)}=\im X_N f_1-\re X_N f_2+\frac{\im(\bar X_Nz^{(N)})}{1+\gamma^{(N)}}Q
$ admitting the same estimate as $X_N$.
Note also that it  follws from our  analysis that for $0\leq \rho\leq 10 t^{-\nu+\varepsilon_1}$,
$0<t\leq T(N)$ ,
\be\label{estimate: 1}
 |\rho^{-l}\partial^k_\rho Z_{in}^{(N)}|\leq C_{k,l}t^{2\nu}<\rho>^{1-l-k}\ln(2+\rho),
\quad k\in \N, \quad  l\leq (3-k)_+.\ee
As a consequence, we obtain the following result.
\begin{lemma}\label{in}
There exists $T(N)>0$ such that for any  $0<t\leq T(N)$ the following holds.\\
(i) The profile  $Z_{in}^{(N)}(\rho,t)$ verifies
\begin{align}
&\|\partial_\rho Z_{in}^{(N)}(t)\|_{L^2(\rho d \rho, 0\leq \rho\leq 10t^{-\nu+\varepsilon_1})}
\leq Ct^\nu,\label{p1.1.11}\\
&\|\rho^{-1}Z_{in}^{(N)}(t)\|_{L^2(\rho d \rho, 0\leq  \rho\leq10t^{-\nu+\varepsilon_1})}
\leq Ct^\nu,\label{p1.1.111}\\
&\| Z_{in}^{(N)}(t)\|_{L^\infty( 0\leq  \rho\leq10t^{-\nu+\varepsilon_1})}+\|\rho\partial_\rho Z_{in}^{(N)}(t)\|_{L^\infty( 0\leq  \rho\leq10t^{-\nu+\varepsilon_1})}
\leq Ct^\nu,\label{p1.1.12}\\
&\|\rho^{-l}\partial_\rho^kZ_{in}^{(N)}(t)\|_{L^2(\rho d \rho, 0\leq \rho\leq10t^{-\nu+\varepsilon_1})}
\leq Ct^{2\nu}(1+|\ln t|),\quad k+l=2,
\label{p1.1.21}\\
&\|\rho^{-l}\partial_\rho^kZ_{in}^{(N)}(t)\|_{L^2(\rho d \rho, 0\leq \rho\leq10t^{-\nu+\varepsilon_1})}
\leq Ct^{2\nu},\quad  k+l\geq 3, \,l\leq(3-k)_+,
\label{p1.1.22}\\
&\|\partial_\rho Z_{in}^{(N)}(t)\|_{L^\infty( 0\leq  \rho\leq10t^{-\nu+\varepsilon_1})}
+\|\rho^{-1}Z_{in}^{(N)}(t)\|_{L^\infty(0\leq  \rho\leq10t^{-\nu+\varepsilon_1})}\leq Ct^{2\nu}(1+|\ln t|),
\label{p1.1.31}\\
&\|\rho^{-l}\partial_\rho^kZ_{in}^{(N)}(t)\|_{L^\infty(0\leq  \rho\leq10t^{-\nu+\varepsilon_1})}\leq C t^{2\nu},\quad 2\leq l+k,\, l\leq (3-k)_+.\label{p1.1.41}
\end{align}
(ii) The error ${\cal R}_{in}^{(N)}$ admits the estimates
\be\label{p1.21}
\begin{split}
&\|\rho^{-l}\partial_\rho^k{\mathcal R}_{in}^{(N)}(t)\|_{L^2(\rho d \rho, 0\leq \rho\leq10t^{-\nu+\varepsilon_1})}\leq t^{N\varepsilon_1},\quad 0\leq l+k\leq 3,\\
&\|\rho^{-l}\partial_\rho^k\partial_t{\mathcal R}_{in}^{(N)}(t)\|_{L^2(\rho d \rho, 0\leq \rho\leq10t^{-\nu+\varepsilon_1})}\leq t^{N\varepsilon_1}, \quad 0\leq k+l\leq 1,
\end{split}
\ee
provided $N>\varepsilon_1^{-1}$.

\end{lemma}
\subsection{Self-similar region  $rt^{-1/2}\lesssim1$}

We next consider  the self-similar  region
$\frac{1}{10}t^{\varepsilon_1}\leq rt^{-1/2}\leq 10t^{-\varepsilon_2}$,
where $0<\varepsilon_2<1/2$ to be fixed later.
In this region we expect the solution to be close to $k$. In this regime it will be convenient
to use the stereographic representation of \eqref{1.1}:
$$(v_1, v_2, v_3)=v\rightarrow w={v_1+iv_2\over 1+v_3} \in \C\cup \{\infty\}.$$
Equation \eqref{1.1} is equivalent to
\be\label{st}
iw_t=-\triangle w+r^{-2}w+G(w,\bar w, w_r),\quad G(w,\bar w, w_r)=
{2\bar w\over 1+|w|^2}(w^2_r-r^{-2}w^2).\ee
Slightly more generally,
if $w(r,t)$ is a solution of
\be
iw_t=-\triangle w+r^{-2}w+G(w,\bar w, w_r)+A,\ee
then $v=\big(\frac{2\re w}{1+|w|^2}, \frac{2\im w}{1+|w|^2},\frac{1-|w|^2}{1+|w|^2})\in S^2$
solves
\be\label{1.11}
v_t=v\times (\triangle v+{R^2\over r^2}v)+{\mathcal A},\ee
with ${\mathcal A}=({\mathcal A}_1,{\mathcal A}_2, {\mathcal A}_3)$  given by
$${\mathcal A}_1+i{\mathcal A}_2=-2i\frac{A+w^2\bar A}{(1+|w|^2)^2},\quad 
{\mathcal A}_3=\frac{4\im(w\bar A)}{(1+|w|^2)^2}.$$
Consider \eqref{st}.
Write $w $ as 
$$w(r,t)=e^{i\a(t)}W(y, t),\quad y=rt^{-1/2}.$$
Then \eqref{st} becomes
\be\label{st1}
itW_t- \a_0W={\cal L}W+G(W,\bar W, W_y),
\ee
where 
$${\cal L}=-\triangle+y^{-2}+i{1\over 2} y\partial_y.$$
Note that as $y\rightarrow 0$,
 \eqref{1.21} gives the following expansion:
\be\label{1.201}
W(y,t)=\sum_{j\geq 0}\sum_{l=0}^{2j+1}\sum_{i\geq -j+l/2}\alpha(j,i,l)t^{\nu(2j+1)}(\ln y -\nu\ln t)^ly^{2i-1},\ee
where the coefficients $\alpha(j,i,l)$ can be expressed explicitly in terms of $c_{j^\prime, l^\prime}^{k,i^\prime}$, $1\leq k\leq j+i$, $ j^\prime\leq i$, $0\leq l^\prime\leq l$.
This suggests the following ansatz for $W$:
\be\label{1.30}
W(y,t)=\sum_{j\geq 0}\sum_{l=0}^{2j+1}t^{\nu(2j+1)}(\ln y -\nu\ln t)^lW_{j,l}(y).\ee
Substituting \eqref{1.30} into \eqref{st1} one gets the following recurrent system for $W_{j,l}$,
$0\leq l\leq 2j+1,\,\,j\geq 0,$:
\begin{align}
&\left\{\begin{array}{l}
({\cal L}-\mu_0)W_{0,1}=0,\\
({\cal L}-\mu_0)W_{0,0}=
-i(1/2+\nu) W_{0,1}+2y^{-1}\partial_yW_{0,1},\end{array}\right.
\label{1.311}\\
&\left\{\begin{array}{l}
({\cal L}-\mu_j)W_{j,2j+1}={\cal G}_{j,2j+1},\\ 
({\cal L}-\mu_j)W_{j,2j}={\cal G}_{j,2j}-i(2j+1)(1/2+\nu)W_{j,2j+1}+
2(2j+1)y^{-1}\partial_yW_{j,2j+1},\\
({\cal L}-\mu_j)W_{j,l}={\cal G}_{j,l}-i(l+1)(1/2+\nu)W_{j,l+1}\\
+
2(l+1)y^{-1}\partial_yW_{j,l+1}+(l+1)(l+2)y^{-2}W_{j,l+2},\quad 0\leq l\leq 2j-1.
\end{array}\right.\label{1.3111}\end{align}
Here 
$\mu_j=-\a_0+i\nu(2j+1)$, 
and ${\cal G}_{j,l}$ is the contribution of the nonlinear term
$ G(W,\bar W, W_r)$, that depends only on 
$W_{i,n}, \, i\leq j-1$:
$$G(W,\bar W, W_r)=
-\sum_{j\geq 1}\sum_{l=0}^{2j+1}t^{(2j+1)\nu}(\ln y -\nu\ln t)^l{\cal G}_{j,l}(y),$$
$${\cal G}_{j,l}(y)={\cal G}_{j,l}(y;W_{i,n},\, 0\leq n\leq 2i+1,\,0\leq i\leq j-1).$$

One has 
\begin{lemma}\label{l2}
Given coefficients
$a_{j}$, $b_j$, $j\geq 0$,
there exists a unique solution of \eqref{1.311}, \eqref{1.3111}, 
$W_{j,l}\in C^\infty(\R_+^*)$,  $0\leq m\leq 2j+1,\,j\geq 0$,
such that as $y\rightarrow 0$, $W_{j,l}$ has the following asymptotic expansion
\be\label{1.33}W_{j,l}(y)=\sum_{i\geq -j+l/2}d^{j,l}_iy^{2i-1},\ee
with \be\label{1.34} 
d^{j,1}_1=a_j,\quad d^{j,0}_1=b_j.\ee
The asymptotic expansion \eqref{1.33} can be differentiated any number of times with respect to $y$.
\end{lemma}

{\it Proof}.
First note that 
equation $({\cal L}-\mu_j)f=0$ has a basis of solutions
$\{e_j^1, e_j^2\}$ such that\\
\noindent (i)
$e^1_j$ is a  $C^\infty$ odd function, $e_j^1(y)=y+O(y^3)$
as $y\rightarrow 0$;\\
\noindent (ii) $e_j^2 \in C^\infty(\R_+^*)$ and admits the representation:
$$e_j^2(y)=y^{-1}+\kappa_je_j^1(y)\ln y+ \tilde e_j^2(y),\quad \kappa_j=
-\frac{i}{4}-\frac{\mu_j}{2},$$
where $\tilde e_j^2$ is   a $C^\infty$ odd function,
$\tilde e_j^2(y)=O(y^3)$ as $y\rightarrow 0$.\\

Consider \eqref{1.311}.
From
$({\cal L}-\mu_0)W_{0,1}=0$ and \eqref{1.33}, \eqref{1.34}, we get
$$W_{0,1}=a_0e_0^1.$$
Consider the equation for $W_{0,0}$:
$$({\cal L}-\mu_0)W_{0,0}=
-i(1/2+\nu) W_{0,1}+2y^{-1}\partial_yW_{0,1}.$$
 The right hand side 
 has the form:
$2a_0y^{-1}$ + a $C^\infty$ odd function. Therefore, the equation has a unique solution
$W_{0,0}^0$ of the form
$$W_{0,0}^0(y)= d_0y^{-1}+\tilde W_{0,0}^0(y),$$
where
$d_0=\frac{a_0}{k_j}$ and $\tilde W_{0,0}^0$ is a $C^\infty$ odd function,
$\tilde W_{0,0}^0(y)=O(y^3)$ as $y\rightarrow 0$.
Together with \eqref{1.33}, \eqref{1.34}, this gives:
$$W_{0,0}=W_{0,0}^0+b_0e_0^1.$$

Consider the case $j\geq 1$. We have
\be\label{1.32}
({\cal L}-\mu_j)W_{j,l}={\cal F}_{j,l},\quad 0\leq  l\leq 2j+1,\ee
where
\begin{equation}\label{1.351}\begin{split}
&{\cal F}_{j,2j+1}={\cal G}_{j,2j+1},\\
&{\cal F}_{j,2j}={\cal G}_{j,2j}-i(2j+1)(1/2+\nu)W_{j,2j+1}+
2(2j+1)y^{-1}\partial_yW_{j,2j+1},\\
&{\cal F}_{j,l}={\cal G}_{j,l}-i(l+1)(1/2+\nu)W_{j,l+1}\\
&+
2(l+1)y^{-1}\partial_yW_{j,l+1}+(l+1)(l+2)y^{-2}W_{j,l+2},\quad 0\leq l\leq 2j-1.
\end{split}\end{equation}
The resolution of \eqref{1.32} is based on the following obvious ODE lemma:
\begin{lemma}\label{l*}
Let $F$ be a $C^\infty({\Bbb R}_+^*)$ function of the form
$$F(y)=\sum_{j= k}^0F_jy^{2j-1}+\tilde F(y),$$
where $ \tilde F$ is a $C^\infty $ odd function and $k\leq -1$.
 Then there exists a unique constant $A$ such that
the equation $({\cal L}-\mu_j)u=F+Ay^{-3}$ has a solution
$u\in C^\infty({\Bbb R}_+^*)$ with the following behavior as 
$y\rightarrow 0$:
$$u(y)=\sum_{j\geq k+1}u_jy^{2j-1},\quad u_1=0.$$
\end{lemma}
More precisely, we proceed as follows.
Assume that $W_{i, n}$, $0\leq n\leq 2i+1$, $i\leq j-1$ has the prescribed behavior \eqref{1.33}, \eqref{1.34}.
Then it is not difficult to check that ${\cal G}_{j,l}$ admit the following expansion as $y\rightarrow 0$:
\be\label{1.36}\begin{split}
&{\cal G}_{j,2j+1}(y)=\sum_{i\geq 1} g_{j,2j+1}^i y^{2i-1},\\
&{\cal G}_{j,2j}(y)=\sum_{i\geq 0} g_{j,2j}^i y^{2i-1},\\
&{\cal G}_{j,l}(y)=\sum_{i\geq -j+l/2-1} g_{j,l}^i y^{2i-1},\quad l\leq 2j-1.
\end{split}
\ee
Consider $W_{j, 2j+1}$. From $({\cal L}-\mu_j)W_{j, 2j+1}={\cal G}_{j,2j+1}$
we get 
\be\label{1.37}
W_{j, 2j+1}=W_{j, 2j+1}^0+ c_0e_j^1,\ee
where $W_{j, 2j+1}^0$ is a unique  $C^\infty$ odd solution of $({\cal L}-\mu_j)f={\cal G}_{j,2j+1}$
that satisfies $W_{j, 2j+1}^0(y)=O(y^3)$ as $y\rightarrow 0$.
The constant $c_0$ remains undetermined at this stage.

Consider ${\cal F}_{j,2j}$.  
It has the form: 
$(g_{j,2j}^0+2(2j+1)c_0)y^{-1}$ + a  $C^\infty$ odd function.
Therefore, for $W_{j,2j}$ we obtain
\be\label{1.38}
W_{j,2j}=W_{j,2j}^0+c_1e_j^1,\ee
where $W_{j,2j}^0$ is a unique solution of $({\cal L}-\mu_j)f={\cal F}_{j,2j}$,
that satisfies as $y\rightarrow 0$,
\be\label{1.391}
W_{j,2j}^0=d_1y^{-1}+O(y^3),\quad d_1=\frac{g_{j,2j}^0+2(2j+1)c_0}{2k_j}.\ee
Similarly to $c_0$, the constant $c_1$ is arbitrary here.

Consider ${\cal F}_{j,2j-1}$. It follows from \eqref{1.351},   \eqref{1.36},  \eqref{1.37},
 \eqref{1.38},  \eqref{1.391}
 that
$${\cal F}_{j,2j-1}=(g_{j,2j-1}^{-1}-4jd_1)y^{-3}+ const\, y^{-1}+{\rm an }\,\,C^\infty \,\, {\rm odd}\,\,{\rm function }.$$
The equation 
$({\cal L}-\mu_j)W_{j,2j-1}={\cal F}_{j,2j-1}$ has a solution of form \eqref{1.33} iff
$$g_{j,2j-1}^{-1}-4jd_1=0,$$
which gives
$$c_0=\frac{k_jg_{j,2j-1}{-1}-2jg_{j,2j}^0}{4j(2j+1)}.$$
With this choice of $c_0$ one gets
$$W_{j,2j-1}=W_{j,2j-1}^0+c_2e_j^1,$$
where $W_{j,2j-1}^0$ is a unique solution of $({\cal L}-\mu_j)f={\cal F}_{j,2j-1}$,
that satisfies as $y\rightarrow 0$,
$$W_{j,2j-1}^0=const \,y^{-1}+O(y^3).$$
Continuing the procedure one successively finds $W_{j, 2j-2},  ...., W_{j, 0}$
in the form $W_{j, 2j+1-k}= W_{j, 2j+1-k}^0+c_ke_j^1$, $k\leq 2j+1$,
where $ W_{j, 2j+1-k}^0$ is an unique solution of $({\cal L}-\mu_j)f={\cal F}_{j,2j+1-k}$,
that as $y\rightarrow 0$ has an asymptotic expansion of the form \eqref{1.33}
with vanishing coefficients $d^{j,l}_1$.
The constant $c_k$, $k\leq 2j-1$, is
 determined uniquely by the solvability condition of the equation
for $W_{j,2j-k-1}$ (see lemma \ref{l*}). Finally,  $c_{2j+1}$,  $c_{2j+2}$
are given by \eqref{1.34}:
$$c_{2j+1}=a_j,\quad c_{2j+2}=b_j.  ~~~~~~~~~~~~~~~~~~~~~~~~~~~~~~~~\square$$
We denote by $W^{ss}_{j,l}(y)$, $0\leq l\leq 2j+1$, $j\geq 0$, the solution of \eqref{1.311},
\eqref{1.3111}
 given by lemma \ref{l2}
with $a_j=\alpha(j,1,1),\,\, b_j=\alpha (j,1,0)$, see \eqref{1.201}.
By uniqueness we have as $y\rightarrow 0$, 
\be\label{1.331}W^{ss}_{j,l}(y)=\sum_{i\geq -j+l/2}\alpha(j,i,l)y^{2i-1},\ee
We next study the behavior of $W_{j,l}^{ss}$,  $0\leq l\leq 2j+1$, $j\geq 0$, at infinity.
One has 
\begin{lemma}\label{l3}
Given coefficients $a_{j,l}$, $b_{j,l}$, $0\leq l\leq 2j+1$, $j\geq 0$, there exists a unique solution of 
 \eqref{1.311}, \eqref{1.3111} of the following form.
\be\label{1.3511}
W_{0,l}=W_{0,l}^0+W_{0,l}^1,\quad l=0,1,\ee
\be\label{1.3512}
W_{j,l}=W_{j,l}^0+W_{j,l}^1+W_{j,l}^2,\quad 0\leq l\leq 2j+1,\,\,j\geq 1,\ee
where $(W_{j,l}^i)_{{0\leq l\leq 2j+1\atop j\geq 1}}$, $i=0,1$, are two solutions of
\eqref{1.311}, \eqref{1.3111} that, as $y\rightarrow \infty$, have the following asymptotic expansion
\be\label{1.3517}\begin{split}
&\sum_{l=0}^{2j+1}(\ln y -\nu\ln t)^l
W_{j,l}^{i,}(y)
=\sum_{l=0}^{2j+1}(\ln y +(-1)^{i}\ln t/2)^l
\hat W_{j,l}^{i}(y),\quad i=0,1,\\
&\hat W_{j,l}^{0}(y)=y^{2i\a_0+2\nu(2j+1)}
\sum_{k\geq 0}\hat w_{k}^{j,l,0}
y^{-2k},\\
&\hat W_{j,l}^{1}(y)=e^{iy^2/4}y^{-2i\a_0-2-2\nu(2j+1)}
\sum_{k\geq 0}\hat w_{k}^{j,l,-1}
y^{-2k},\end{split}\ee
with
\be\label{condition}
\hat w_{0}^{j,l,0}=a_{j,l},\quad \hat w_{0}^{j,l,-1}=b_{j,l}.
\ee
Finally, the interaction part $ W_{j,l}^2$ can be written as
\be\label{1.355}W_{j,l}^2(y)=\sum_{-j-1\leq m\leq j}e^{-imy^2/4}y^{2i\a_0(2m+1)}W_{j,l,m}(y),
\end{equation}
where $W_{j,l,m}$ have the following asymptotic expansion
as $y\rightarrow \infty$:
\be\label{1.3513}\begin{split}
&W_{j,l,m}(y)=\sum_{k\geq m+2}\sum_{{m-j\leq i\leq j-m\atop j-m-i\in 2{\Bbb Z}}}\sum_{s=0}^{2j+1-l}w_{k,i,s}^{j,l,m}
y^{2\nu(2i+1)-2k}(\ln y)^s,\quad m\geq 1,\\
&W_{j,l,m}(y)=\sum_{k\geq -m}\,\sum_{{-j-m-2\leq i\leq j+m\atop j-m-i\in 2{\Bbb Z}}}\sum_{s=0}^{2j+1-l}w_{k,i,s}^{j,l,m}
y^{2\nu(2i+1)-2k}(\ln y)^s,\quad m\leq -2,\\
&W_{j,l,0}(y)=
\sum_{k\geq 1}\sum_{{-j\leq i\leq j-2\atop  j-i\in 2{\Bbb Z}}}\sum_{s=0}^{2j+1-l}w_{k,i,s}^{j,l,0}
y^{2\nu(2i+1)-2k}(\ln y)^s,\\
&W_{j,l, -1}(y)=
\sum_{k\geq 1}\sum_{{-j+1\leq i\leq j-1\atop j-i+1\in 2{\Bbb Z}} }\sum_{s=0}^{2j+1-l}w_{k,i,s}^{j,l,-1}
y^{2\nu(2i+1)-2k}(\ln y)^s.\end{split}\ee
The asymptotic expansion \eqref{1.3517}, \eqref{1.3513} can be differentiated any number of times with respect to $y$.

Any solution of  \eqref{1.311}, \eqref{1.3111} has form \eqref{1.3511},  \eqref{1.3512},
 \eqref{1.3517},  \eqref{1.355}, \eqref{1.3513}.
\end{lemma}
\par {\it Proof}.
First note that 
equation $({\cal L}-\mu_j)f=0$ has a basis of solutions
$\{f_j^1, f_j^2\}$  with the following behavior at infinity:
$$f_j^1(y)=y^{2i\a_0+2\nu(2j+1)}\sum_{k\geq 0}f^{k}_{j,1}y^{-2k}, \quad
f_j^2(y)=e^{iy^2/4}y^{-2i\a_0-2\nu(2j+1)-2}\sum_{k\geq 0}f^{k}_{j,2}y^{-2k},$$
$ f^{0}_{j,1}= f^{0}_{j,2}=1$.
As a consequence, the homogeneous system
\be\label{1.3710}\begin{split}
&({\cal L}-\mu_j)g_{2j+1}=0,\\ 
&({\cal L}-\mu_j)g_{2j}=-i(2j+1)(1/2+\nu)g_{2j+1}+
2(2j+1)y^{-1}\partial_yg_{2j+1},\\
&({\cal L}-\mu_j)g_{l}=-i(l+1)(1/2+\nu)g_{l+1}\\
&+
2(l+1)y^{-1}\partial_yg_{l+1}+(l+1)(l+2)y^{-2}g_{l+2},\quad 0\leq l\leq 2j-1.
\end{split}\ee
has a basis of solutions $\{{\bold g}_{j}^{i,m}\}_{{i=1,2, \atop m=0,\dots , 2j+1}}$,
$${\bold g}_{j}^{i,m}=( g_{j,0}^{i,m}, \dots,  g_{j,2j+1}^{i,m}),\quad  0\leq m\leq 2j+1,\,\,i=1,2,$$
defined by
\be\label{1.49}
\sum_{l=0}^{2j+1}(\ln y -\nu\ln t)^lg^{i,m}_{j, l}(y)=
\sum_{l=0}^{2j+1}(\ln y +(-1)^{i-1}\ln t/2)^l\xi^{i,m}_{j,l}(y),\ee
where $(\xi^{i,m}_{j,l})_{l=0,\dots, 2j+1}$ is the unique solution of
\be\label{1.50}\begin{split}
&({\cal L}-\mu_j)\xi_{2j+1}=0,\\ 
&({\cal L}-\mu_j)\xi_{2j}=-i(2j+1)(i-1)\xi_{2j+1}+
2(2j+1)y^{-1}\partial_y\xi_{2j+1},\\
&({\cal L}-\mu_j)\xi_{l}=-i(l+1)(i-1)\xi_{l+1}+
2(l+1)y^{-1}\partial_y\xi_{l+1}\\
&+(l+1)(l+2)y^{-2}\xi_{l+2},\quad 0\leq l\leq 2j-1,
\end{split}\ee
verifying
\be\label{1.51}\begin{split}
&\xi_{j,l}^{i,m}(y)=0, \quad l>2j+1-m,\\
&\xi_{j,2j+1-m}^{i,m}(y)=f_j^i(y),\\
&\xi_{j,l}^{1,m}(y)=y^{2i\a_0+2\nu(2j+1)}
\sum_{k\geq 2j+1-l-m}\xi_{l,k}^1y^{-2k},\quad y\rightarrow +\infty,\\
&\xi_{j,l}^{2,m}(y)=e^{iy^2/4}y^{-2i\a_0-2\nu(2j+1)-2}
\sum_{k\geq 2j+1-l-m}\xi_{l,k}^2y^{-2k}\quad y\rightarrow +\infty.\end{split}\ee

Consider $W_{0, l}$, $l=0,1$.
We have
$$({\cal L}-\mu_0)W_{0,1}=0,$$ 
$$({\cal L}-\mu_0)W_{0,0}=
-i(1/2+\nu) W_{0,1}+2y^{-1}\partial_yW_{0,1},$$
which gives
$$W_{0,l}(y)=\sum_{i=1,2,\atop m=0,1}A_{i,m}g^{i,m}_{0,l} (y),\quad l=0,1,$$
with some constants $A_{i,m}$, $i=1,2$, $m=0, 1$.
It follows from\eqref{1.49},  \eqref{1.51} that
 $W_{0, l}$, $l=0,1$ have the form \eqref{1.3511}, \eqref{1.3517}
with
$\hat w^{j,l,0}_{0}=A_{1, 1-l}$,  $\hat w^{j,l,-1}_{0}=A_{2, 1-l}$, $l=0,1$,
which together with \eqref{condition} gives
$A_{1,m}=a_{0,1-m}$, $ A_{2,m}=b_{0, 1-m}$, $m=0,1$.

We next consider $j\geq 1$.
Assume that $W_{i, n}$, $0\leq n\leq 2i+1$, $i\leq j-1$ has the prescribed behavior \eqref{1.3512}, \eqref{1.3517},\eqref{1.355}, \eqref{1.3513}.
Then it is not difficult to check that ${\cal G}_{j,l}$ 
has the form
\be\label{1.40}
{\cal G}_{j,l}(y)=\sum_{-j-1\leq m\leq j}e^{-imy^2/4}y^{2i\a_0(2m+1)}{\cal G}_{j,l}^m(y),
\ee
where ${\cal G}_{j,l}^m$, $m=0,-1$, are given by
\be\label{1.400}
\begin{split}
&{\cal G}_{j,l}^m(y)={\cal G}_{j,l}^{m,0}(y)+{\cal G}_{j,l}^{m,1}(y),\quad m=0,-1,\\
&{\cal G}_{j,l}^{0,0}(y)={\cal G}_{j,l}(y;W^0_{i,n}, 0\leq n\leq 2i+1, 0\leq i\leq j-1),\\
&e^{iy^2/4}{\cal G}_{j,l}^{-1,0}(y)={\cal G}_{j,l}(y;W^1_{i,n}, 0\leq n\leq 2i+1, 0\leq i\leq j-1),
\end{split}\ee
and admit  the following asymptotic expansions
and as $y\rightarrow\infty$:
\be\label{1.402}
\begin{split}
&{\cal G}_{j,l}^{0,0}(y)=\sum_{k\geq 1}\sum_{s=0}^{2j+1-l}T_{k,j,s}^{j,l,0}
y^{2\nu(2j+1)-2k}(\ln y)^s,\\
&{\cal G}_{j,l}^{0,1}(y)=
\sum_{k\geq 2}\sum_{{-j\leq i\leq j-2\atop  j-i\in 2{\Bbb Z}}}\sum_{s=0}^{2j+1-l}T_{k,i,s}^{j,l,0}
y^{2\nu(2i+1)-2k}(\ln y)^s,
\end{split}\ee
\be\label{1.403}
\begin{split}
&{\cal G}_{j,l}^{-1,0}(y)=\sum_{k\geq 2}\sum_{s=0}^{2j+1-l}T_{k,-j-1,s}^{j,l,-1}
y^{-2\nu(2j+1)-2k}(\ln y)^s,\\
&{\cal G}_{j,l}^{-1,1}(y)
\sum_{k\geq 1}\sum_{{-j+1\leq i\leq j-1\atop  j-i+1\in 2{\Bbb Z}}}\sum_{s=0}^{2j+1-l}T_{k,i,s}^{j,l,-1}
y^{2\nu(2i+1)-2k}(\ln y)^s.\end{split}\ee
Finally,  ${\cal G}_{j,l}^m$, $m\neq 0,-1$,  have the following behavior as $y\rightarrow \infty$
\be\label{1.401}
\begin{split}
&{\cal G}_{j,l}^m(y)=
\sum_{k\geq m+1}\sum_{{m-j\leq i\leq j-m\atop j-m-i \in 2{\Bbb Z}}}\sum_{s=0}^{2j+1-l}T_{k,i,s}^{j,l,m}
y^{2\nu(2i+1)-2k}(\ln y)^s,\quad m\geq 1,\\
&{\cal G}_{j,l}^m(y)=\sum_{k\geq |m|-1}\sum_{{-j-m-2\leq i\leq j+m\atop j+m-i \in{\Bbb Z}}}\sum_{s=0}^{2j+1-l}T_{k,i,s}^{j,l,m}
y^{2\nu(2i+1)-2k}(\ln y)^s,\quad m\leq -2,
\end{split}\ee
Therefore, integrating  \eqref{1.3111}, one gets
\begin{equation}\begin{split}\label{1.335}
&W_{j,l}=\tilde W_{j,l}+\sum_{{i=1,2\atop m=0,\dots, 2j+1}}A_{i,m}{\bold g}_{j,l}^{i,m},\\
&\tilde W_{j,l}(y)=\sum_{-j-1\leq m\leq j}e^{-imy^2/4}y^{2i\a_0(2m+1)}\tilde W_{j,l}^m(y),
\end{split}
\end{equation}
where $e^{-imy^2/4}y^{2i\a_0(2m+1)}\tilde W_{j,l}^m$ is a unique solution of  \eqref{1.3111} 
with ${\mathcal G}_{j,l}$ replaced by $e^{-imy^2/4}y^{2i\a_0(2m+1)}{\cal G}_{j,l}^m(y)$,
that has 
the following behavior  as $y\rightarrow +\infty$:
\be\label{1.35}\begin{split}
&\tilde W_{j,l}^m(y)=\sum_{k\geq m+2}\,\sum_{{m-j\leq i\leq j-m\atop j-m-i\in 2{\Bbb Z}}}\sum_{s=0}^{2j+1-l}\tilde w_{k,i,s}^{j,l,m}
y^{2\nu(2i+1)-2k}(\ln y)^s,\quad m\geq 1,\\
&\tilde W_{j,l}^m(y)=\sum_{k\geq -m}\,\sum_{{-j-m-2\leq i\leq j+m\atop j-m-i\in 2{\Bbb Z}}}\sum_{s=0}^{2j+1-l}\tilde w_{k,i,s}^{j,l,m}
y^{2\nu(2i+1)-2k}(\ln y)^s,\quad m\leq -2.\\
\end{split}
\end{equation}
Finally for $m=0,-1$ one has:
\be\label{1.35*}\begin{split}
&\tilde W_{j,l}^0(y)=\tilde W_{j,l}^{0,0}(y)+\tilde W_{j,l}^{0,1}(y),\\
&\tilde W_{j,l}^{-1}(y)=\tilde W_{j,l}^{-1,0}(y)+\tilde W_{j,l}^{-1,1}(y),
\end{split}
\ee
where
$\tilde W_{j,l}^{0,i}$ and $e^{iy^2/4}y^{-2i\alpha_0}\tilde W_{j,l}^{-1,i}$ are solutions 
of
 \eqref{1.3111} with ${\cal G}_{j,l}$ replaced by ${\cal G}_{j,l}^{0,i}$
and $y^{-2i\alpha_0}e^{iy^2/4}{\cal G}_{j,l}^{-1,i}$  respectively,  with the following asymptotics as $y\rightarrow \infty$:
\be\label{1.535}\begin{split}
&\tilde W_{j,l}^{0,0}(y)=
\sum_{k\geq 1}\sum_{s=0}^{2j+1-l}\tilde w_{k,j,s}^{j,l,0}
y^{2\nu(2j+1)-2k}(\ln y)^s,\\
&\tilde W_{j,l}^{0,1}(y)=
\sum_{k\geq 1}\sum_{{-j\leq i\leq j-2\atop  j-i\in 2{\Bbb Z}}}\sum_{s=0}^{2j+1-l}\tilde w_{k,i,s}^{j,l,0}
y^{2\nu(2i+1)-2k}(\ln y)^s,\\
&\tilde W_{j,l}^{-1,0}(y)=
\sum_{k\geq 2}\sum_{s=0}^{2j+1-l}\tilde w_{k,-j-1,s}^{j,l,-1}
y^{-2\nu(2j+1)-2k}(\ln y)^s,\\
&\tilde W_{j,l}^{-1,1}(y)=
\sum_{k\geq 1}\sum_{{-j+1\leq i\leq j-1\atop  j-i\in 2{\Bbb Z}}}\sum_{s=0}^{2j+1-l}
\tilde w_{k,i,s}^{j,l,-1}
y^{2\nu(2i+1)-2k}(\ln y)^s.\end{split}\ee
Clearly,
$W_{j,l}^0 =\tilde W_{j,l}^{0,0}+\sum_{m=0}^{2j+1}A_{1,m}{\bold g}_{j,l}^{1,m}$,
and $W_{j,l}^1 =e^{-imy^2/4}\tilde W_{j,l}^{-1,0}+
\sum_{m=0}^{2j+1}A_{2,m}{\bold g}_{j,l}^{2,m}$ are solutions of 
 \eqref{1.3111} with ${\cal G}_{j,l}$ replaced by ${\cal G}_{j,l}^{0,0}={\cal G}_{j,l}(W^0_{i,n}, \, i\leq j-1)$ and $e^{iy^2/4}{\cal G}_{j,l}^{-1,0}={\cal G}_{j,l}(W^1_{i,n}, \, i\leq j-1)$ respectively.
As a consequence, $W_{j,l}^i$, $i=0,1$, $0\leq l\leq 2j+1$, have the form \eqref{1.3517} 
with $\hat w_0^{j,l ,i}=A^{1-i}_{2j+1-l}$, $i=0, -1$, $l=0,\dots, 2j+1$,
which together with \eqref{condition} gives
$A_{1,m}=a_{j,2j+1-m}$, $ A_{2,m}=b_{j, 2j+1-m}$, $m=0,\dots, 2j+1$.
 $~~~~~~~~~~~~~~~~~~~\square$

Let $W_{in}^{(N)}(y,t)$ be the  
 the stereographic representation of 
$V_{in}^{(N)}(t^{-\nu}y,t)=(V_{in,1}^{(N)}(t^{-\nu}y,t), V_{in,2}^{(N)}(t^{-\nu}y,t), V_{in,3}^{(N)}(t^{-\nu}y,t))$:
$$ W_{in}^{(N)}(y,t)=\frac{V_{in,1}^{(N)}(t^{-\nu}y,t)+i V_{in,2}^{(N)}(t^{-\nu}y,t)}{ 1+V_{in,3}^{(N)}(t^{-\nu}y,t)} .$$
For $N\geq 2$ define
\begin{equation*}
\begin{split}
&W_{ss}^{(N)}(y,t)=\sum_{j=0}^N\sum_{l=0}^{2j+1}t^{\nu(2j+1)}(\ln \rho)^lW^{ss}_{j,l}(y),\\
&A_{ss}^{(N)}=-it\partial_tW_{ss}^{(N)}+ \a_0W_{ss}^{(N)}+{\cal L}W_{ss}^{(N)}+
G(W_{ss}^{(N)},\bar W_{ss}^{(N)}, \partial_y W_{ss}^{(N)}),\\
&V_{ss}^{(N)}(\rho,t)=\big(\frac{2\re  W_{ss}^{(N)}}{1+| W_{ss}^{(N)}|^2}, \frac{2\im  W_{ss}^{(N)}}{1+| W_{ss}^{(N)}|^2},
\frac{1-| W_{ss}^{(N)}|^2}{1+| W_{ss}^{(N)}|^2}), \quad \rho=t^{-\nu}y,\\
&Z_{ss}^{(N)}(\rho,t)=V^{(N)}_{ss}(\rho,t)-Q(\rho).
\end{split}
\end{equation*}
Fix $\varepsilon_1=\frac{\nu}{2}$.  Then, as a direct consequence of the previous analysis, 
we obtain the following result.
\begin{lemma}\label{ss}
For $0<t\leq T(N)$ the following holds.\\
(i) For any $k,\, l$,and  $ \frac{1}{10}t^{\varepsilon_1}\leq  y \leq {10}t^{\varepsilon_1}$,
one has
\be\label{ss.1}
|y^{-l}\partial_y^{k}\partial_t^i(W_{ss}^{(N)}-W_{in}^{(N)})|\leq C_{k,l,i}t^{\nu(N+1-\frac{l+k}{2})-i},\quad i=0,1.
\ee
(ii)The profile $Z_{ss}^{(N)}$ verifies
\begin{align}
&\|\partial_\rho Z_{ss}^{(N)}(t)(\|_{L^2(\rho d \rho,  \frac{1}{10}t^{-\nu+\varepsilon_1}\leq \rho\leq 10t^{-\nu-\varepsilon_2})}
\leq Ct^{\eta},\label{ss.1.11}\\
&\|\rho^{-1}Z_{ss}^{(N)}(t)\|_{L^2(\rho d \rho,  \frac{1}{10}t^{-\nu+\varepsilon_1}\leq \rho\leq 10t^{-\nu-\varepsilon_2})}
\leq Ct^{\eta},\label{p1.1.111111}\\
&\|Z_{ss}^{(N)}(t)\|_{L^\infty( \frac{1}{10}t^{-\nu+\varepsilon_1}\leq \rho\leq 10t^{-\nu-\varepsilon_2})}\leq Ct^{\eta},\label{ss.1.12}\\
&\|\rho\partial_\rho Z_{ss}^{(N)}(t)\|_{L^\infty( \frac{1}{10}t^{-\nu+\varepsilon_1}\leq \rho\leq 10t^{-\nu-\varepsilon_2})}
\leq Ct^{\eta},\label{p1.1.1212}\\
&\|\rho^{-l}\partial_\rho^kZ_{ss}^{(N)}(t)\|_{L^2(\rho d \rho, \frac{1}{10}t^{-\nu+\varepsilon_1}\leq \rho\leq 10t^{-\nu-\varepsilon_2})}
\leq Ct^{\nu+\frac12+\eta},\quad k+l=2,
\label{p1.1.2112}\\
&\|\rho^{-l}\partial_\rho^kZ_{in}^{(N)}(t)\|_{L^2(\rho d \rho, \frac{1}{10}t^{-\nu+\varepsilon_1}\leq \rho\leq 10t^{-\nu-\varepsilon_2})}
\leq Ct^{2\nu},\quad  k+l\geq 3,
\label{p1.1.2222}\\
&\|\rho^{-l}\partial_\rho^kZ_{ss}^{(N)}(t)\|_{L^\infty(\frac{1}{10}t^{-\nu+\varepsilon_1}\leq \rho\leq 10t^{-\nu-\varepsilon_2})} 
\leq C
t^{\nu+\eta}, \quad k+l=1,
\label{p1.1.3131}\\
&\|\rho^{-l}\partial_\rho^kZ_{ss}^{(N)}(t)\|_{L^\infty( \frac{1}{10}t^{-\nu+\varepsilon_1}\leq \rho\leq 10t^{-\nu-\varepsilon_2})}\leq C t^{2\nu},\quad 2\leq l+k.\label{ss1.1.41}
\end{align}
Here and below $\eta$ stands for small positive constants depending on $\nu$ and $\varepsilon_2$,
that may change from line to line.\\
(iii) The error $A_{ss}^{(N)}$ admits the estimate
\be\label{p1.1.2121}
\|y^{-l}\partial_y^k\partial_t^iA_{ss}^{(N)}(t)\|_{L^2(yd y, \frac{1}{10}t^{\varepsilon_1}\leq y\leq10t^{-\varepsilon_2})}\leq t^{\nu N(1-2\varepsilon_2)-i},
\quad 0\leq l+k\leq 4,\quad i=0,1.
\ee
\end{lemma}

\subsection{Remote region  $r\sim 1$}
We next consider  the remote  region
$t^{-\varepsilon_2}\leq rt^{-1/2}$.
Consider the formal solution $\sum_{j\geq 0}\sum_{l=0}^{2j+1}t^{\nu(2j+1)}(\ln y -\nu\ln t)^lW^{ss}_{j,l}(y)$ constructed in the previous subsection.
By lemma \ref{l3}, it has form \eqref{1.3511}, \eqref{1.3512}, \eqref{1.3517}, \eqref{1.355}, \eqref{1.3513},
with some coefficients $\hat w_k^{j,l,i}$, $w_{k,i,s}^{j,l,m}$. Note that in the limit $y\rightarrow \infty$, $r\rightarrow 0$, the main order terms of the expansion
$\sum_{j\geq 0}\sum_{l=0}^{2j+1}t^{\nu(2j+1)+i\alpha_0}(\ln y -\nu\ln t)^lW^{ss}_{j,l}(t^{-1/2}r)$ are given by
\begin{equation}\label{remote1}\begin{split}
&\sum_{j\geq 0}\sum_{l=0}^{2j+1}t^{\nu(2j+1)+i\alpha_0}(\ln y -\nu\ln t)^lW^{ss}_{j,l}(t^{-1/2}r)
\sim \\
&\sum_{k\geq 0}{t^k\over r^{2k}}\sum_{j\geq 0}\sum_{l=0}^{2j+1}\hat w_k^{j,l,0}
(\ln r)^lr^{2i\alpha_0+2\nu(2j+1)}\\
&+
{e^{{ir^2\over 4t}}\over t}\sum_{k\geq 0}{t^k\over r^{2k}}\sum_{j\geq 0}
\sum_{l=0}^{2j+1}\hat w_k^{j,l,0}\left({r\over t}\right)^{-2i\a_0-2\nu(2j+1)-2}\left(\ln \left({r\over t}\right)\right)^l,
\end{split}
\end{equation}
which means that in region $t^{-\varepsilon_2}\leq rt^{-1/2}$
we have to look for the solution of \eqref{st}
as a perturbation
of the time independent profile
$$\sum_{j\geq 0}\sum_{l=0}^{2j+1}\beta_0(j,l)
(\ln r)^lr^{2\nu(2j+1)},$$
 with $\beta_0(j,l)=\hat w_0^{j,l,0}$.

Let $\theta\in C_0^\infty({\Bbb R})$, 
$\theta(\xi)=\left\{\begin{array}{cr}& 1, \quad |\xi|\leq 1,\\&0,\quad |\xi|\geq 2.\end{array}\right.$
For $N\geq 2$, and $\delta>0$  we define
$$f_0(r)\equiv f_0^{(N)}(r)=\theta(\delta^{-1}r)\sum_{j= 0}^N\sum_{l=0}^{2j+1}\beta_0(j,l)
(\ln r)^lr^{2i\alpha_0+2\nu(2j+1)}.$$
Note that $e^{i\theta}f_0\in H^{1+2\nu-}$ and 
\be\label{f_0}
\|e^{i\theta}f_0\|_{\dot H^s}\leq C\delta^{1+2\nu-s},\quad \forall \, 0\leq s<1+2\nu.
\ee

Write $w(r,t)=f_0(r)+\chi(r,t)$. Then $\chi$ solves
\be\label{r1}\begin{split}
&i\chi_t=-\triangle \chi+r^{-2}\chi+{\cal V}_0\partial_r\chi+{\cal V}_1\chi+{\cal V}_2\bar\chi+{\mathcal N}+{\cal D}_0,\\
&{\cal V}_0=\frac{4\bar f_0\partial_rf_0}{1+|f_0|^2},\quad {\cal V}_1=
-\frac{2|f_0|^2(2+|f_0|^2)}{r^2(1+|f_0|^2)^2}-
\frac{2\bar f_0^2(\partial_rf_0)^2}{(1+|f_0|^2)^2},\\
& {\cal V}_2=\frac{2(r^2(\partial_rf_0)^2-f_0^2)}{r^2(1+|f_0|^2)^2},\\
&{\cal D}_0=(-\triangle +r^{-2})f_0+G(f_0, \bar f_0,\partial_r f_0).
\end{split}\end{equation}
Finally,  ${\mathcal N}$ contains the terms that are at least quadratic in $\chi$ and  it has the form
\be\label{r2}\begin{split}
&{\mathcal N}=N_0(\chi,\bar\chi)+\chi_rN_1(\chi,\bar\chi)+\chi_r^2N_2(\chi,\bar\chi),\\
&N_0(\chi,\bar\chi)=G(f_0+\chi, \bar f_0+\bar \chi, \partial_r f_0)-
G(f_0, \bar f_0, \partial_r f_0)-{\cal V}_1\chi-{\cal V}_2\bar\chi,\\
&N_1(\chi,\bar\chi)=\frac{4\partial_rf_0(\bar f_0+\bar\chi)}{1+|f_0+\chi|^2}-{\cal V}_0,\\
&N_2(\chi,\bar\chi)=\frac{2(\bar f_0+\bar\chi)}{1+|f_0+\chi|^2}.
\end{split}
\ee
Accordingly to \eqref{1.3511}, \eqref{1.3512}, \eqref{1.3517}, \eqref{1.355}, \eqref{1.3513}, we look for $\chi $ as 
\be\label{r3}
\chi(r,t)=\sum_{{q\geq 0\atop k\geq 1}}t^{2\nu q +k}\sum_{{-\min\{k,q\}\leq m\leq \min\{(k-2)_+,q\}
\atop q-m\in 2{\Bbb Z}}}\sum_{s=0}^q e^{-im\Phi}(\ln r-\ln t)^sg_{k,q, m,s}(r),\ee
where 
$$\Phi=\frac{r^2}{4t}+2\a_0\ln t+\varphi(r),$$
with $\varphi$ to be chosen later.

Substituting this ansatz to the expressions
$-i\chi_t-\triangle \chi+r^{-2}\chi+{\cal V}_0\partial_r\chi+{\cal V}_1\chi+{\cal V}_2\bar\chi$,
$N$, we get
\begin{equation*}\begin{split}
&-i\chi_t+\triangle \chi-r^{-2}\chi+{\cal V}_0\partial_r\chi+{\cal V}_1\chi+{\cal V}_2\bar\chi\\
&=\sum_{{q\geq 0\atop k\geq 2}}t^{2\nu q +k-2}\sum_{{-\min\{k,q\}\leq m\leq \min\{(k-2)_+,q\}
\atop q-m\in 2{\Bbb Z}}}\sum_{s=0}^q e^{-im\Phi}(\ln r-\ln t)^s\Psi^{lin}_{k,q, m,s},\\
&N_0(\chi,\bar\chi)=\sum_{{q\geq 0\atop k\geq 4}}t^{2\nu q +k-2}\sum_{{-\min\{k,q\}\leq m\leq \min\{(k-2)_+,q\}
\atop q-m\in 2{\Bbb Z}}}\sum_{s=0}^q e^{-im\Phi}(\ln r-\ln t)^s\Psi^{nl,0}_{k,q, m,s},\\
&\chi_rN_1(\chi,\bar\chi)=\sum_{{q\geq 0\atop k\geq 3}}t^{2\nu q +k-2}\sum_{{-\min\{k,q\}\leq m\leq \min\{(k-2)_+,q\}
\atop q-m\in 2{\Bbb Z}}}\sum_{s=0}^q e^{-im\Phi}(\ln r-\ln t)^s\Psi^{nl,1}_{k,q, m,s},\\
&(\chi_r)^2N_2(\chi,\bar\chi)=\sum_{{q\geq 0\atop k\geq 2 }}t^{2\nu q +k-2}\sum_{{-\min\{k,q\}\leq m\leq \min\{(k-2)_+,q\}
\atop q-m\in 2{\Bbb Z}}}\sum_{s=0}^q e^{-im\Phi}(\ln r-\ln t)^s\Psi^{nl,2}_{k,q, m,s},\\
\end{split}\end{equation*}
Here
\be\label{remote:01}
\Psi^{lin}_{k,q, m,s}={m(m+1)r^2\over 4}g_{k,q, m,s}+\Psi^{lin,1}_{k,q, m,s}+\Psi^{lin,2}_{k,q, m,s},\ee
with $\Psi^{lin,1}_{k,q, m,s}$ and  $\Psi^{lin,2}_{k,q, m,s}$
depending respectively on $g_{k-1,q, m,s^\prime}, s^\prime=s, s+1$
and $g_{k-2,q, m,s^\prime}$,  $s^\prime=s, s+1, s+2$ only:
\be\label{remote:10}
\begin{split}
\Psi^{lin,1}_{k,q, m,s}=&-i(2\nu q+k-1-m-2im\a_0)g_{k-1,q, m,s}+i(m+1)(s+1)g_{k-1,q, m,s+1}\\
&+imr(\partial_r-im\varphi^\prime(r)-\frac12{\cal V}_0(r))g_{k-1,q, m,s},
\end{split}\ee
\be\label{remote:100}\begin{split}
&\Psi^{lin,2}_{k,q, m,s}=-e^{im\varphi}\triangle (e^{-im\varphi}g_{k-2,q, m,s})-\frac{2(s+1)}{r}e^{im\varphi}\partial_r(e^{-im\varphi}g_{k-2,q, m,s+1})\\
&-\frac{(s+1)(s+2)}{r^2}g_{k-2,q, m,s+2}
+
{\cal V}_0e^{im\varphi}\partial_r (e^{-im\varphi}g_{k-2,q, m,s})\\&+
(r^{-2}+{\cal V}_1)g_{k-2,q, m,s}+{\cal V}_2\bar g_{k-2,q, -m,s}.
\end{split}\ee
Here and below we use the convention
$g_{k,q, m,s}=0$ if $
(k,q,m,s)\notin \Omega$,
where 
$$\Omega=\{
k\geq 1, q\geq 0, 0\leq s\leq q, q-m\in 2{\Bbb Z}, -\min\{k,q\}\leq m\leq \min\{k-1,q\}\}.$$
The nonlinear terms $\Psi^{nl,i}_{k,q, m,s}$, $i=0,1$, depend only on 
$g_{k^\prime,q^\prime, m^\prime,s^\prime}$ with $k^\prime\leq k-2$. More precisely,
$$\Psi^{nl,0}_{k,q, m,s}=\Psi^{nl,0}_{k,q, m,s}(r;g_{k^\prime,q^\prime, m^\prime,s^\prime},\, k^\prime\leq k-3),$$
$$\Psi^{nl,1}_{k,q, m,s}=\Psi^{nl,1}_{k,q, m,s}(r;g_{k^\prime,q^\prime, m^\prime,s^\prime},\, k^\prime\leq k-2).$$
Finally, $\Psi^{nl,2}_{k,q, m,s}$ has the following structure
\begin{equation}\label{N2}\begin{split}
&\Psi^{nl,2}_{2,q, m,s}=-\delta_{m, -2}\frac{r^2\bar f_0}{2(1+|f_0|^2)}\sum\limits_{{q_1+q_2=q\atop s_1+s_2=s}}g_{1,q_1, -1,s_1}g_{1,q_2, -1,s_2},\\
&\Psi^{nl,2}_{k,q, m,s}=\Psi^{nl,2,0}_{k,q, m,s}+\tilde \Psi^{nl,2}_{k,q, m,s},\quad k\geq 3,\\
&\Psi^{nl,2,0}_{k,q, m,s}=\frac{(m+1)r^2\bar f_0}{1+|f_0|^2}\sum\limits_{{q_1+q_2=q\atop s_1+s_2=s}}g_{1,q_1, -1,s_1}g_{k-1,q_2, m+1,s_2}, 
\end{split}\ee
with
$\tilde \Psi^{nl,2}_{k,q, m,s}$ depending on $g_{k^\prime,q^\prime, m^\prime,s^\prime}$, $k^\prime\leq k-2$
only:
$$ \tilde\Psi^{nl,2}_{k,q, m,s}=
(r;g_{k^\prime,q^\prime, m^\prime,s^\prime},\, k^\prime\leq k-2).$$

Note that 
$$\Psi^{nl,2,0}_{k,q, -1,s}=0,\quad \forall k,q,s.$$
Equation \eqref{r1} is equivalent to
\be\label{R1}
\left\{
\begin{array}{ll}
&\Psi^{lin}_{2,0, 0,0}+{\cal D}_0=0,\\
&\Psi^{lin}_{k,q, m,s}+\Psi^{nl}_{k,q, m,s}=0,\quad   (k,q, m,s)\in \Omega,\,\,(k,q, m,s)\neq (2,0, 0,0), 
\end{array}\right.\ee
Here  $\Psi^{nl}_{k,q, m,s}=\Psi^{nl,0}_{k,q, m,s}+\Psi^{nl,1}_{k,q, m,s}+\Psi^{nl,2}_{k,q, m,s}$.

We view \eqref{R1}
as a recurrent system with respect to $k\geq 1$ of the form
\be\label{R3}
\left\{
\begin{array}{ll}
&\Psi^{lin}_{2,0, 0,0}+{\cal D}_0=0,\\
&\Psi^{lin}_{2,2j, 0,s}=0,\quad (j,s)\neq (0,0),\\
&\Psi^{lin}_{2,2j+1, 1,s}=0,
\end{array}\right.
\ee
and
\be\label{R2}\left\{
\begin{array}{ll}
&\Psi^{lin}_{k+1,q, m,s}+\Psi^{nl}_{k+1,q, m,s}=0, \quad m=0,1\\
&\Psi^{lin}_{k,q, m,s}+\Psi^{nl}_{k,q, m,s}=0, \quad m\neq 0,1
\end{array},\right.\quad k\geq 2.\ee
Consider \eqref{R3}. Choosing $\varphi$ as
\be\label{phi}
\varphi(r)=-i\int_0^r ds \frac{\bar f_0(s)\partial_sf_0(s)- f_0(s)\partial_s\bar f_0(s)}{1+|f_0(s)|^2},
\ee
we can rewrite \eqref{R3} in the following form
\be\label{R4}\left\{
\begin{array}{ll}
&(4\nu j+1)g_{1,2j,0,s}-(s+1)g_{1,2j,0,s+1}=0,\quad (j,s)\neq(0,0),\\
&g_{1,0,0,0}=-iD_0,\\
&r\partial_rg_{1,2j+1,-1,s}+(2\nu (2j+1)+2 +2i\alpha_0-r(\ln(1+|f_0|^2))^\prime)g_{1,2j+1,-1,s}=0.
\end{array}\right.
\ee
Accordingly to  \eqref{remote1}, we solve this system as follows:
\be\label{k=1}\begin{split}
&g_{1,2j,0,s}=0,\quad (j,s)\neq(0,0),\\
&g_{1,0,0,0}=-iD_0,\\
&g_{1, 2j+1, -1, s}=\beta_1(j,s)(1+|f_0|^2)r^{-2i\a_0-2\nu(2j+1)-2},\,\,\, 0\leq s\leq 2j+1,\,\, 0\leq j\leq N,\\
&g_{1,2j+1, -1, s}=0, \quad j>N,
\end{split}\ee
where $\beta_1(j,s)= \hat w_{0}^{j,l,-1}$.

Consider  \eqref{R2}.  We will solve it with the   "zero boundary conditions" at
zero. To formulate the result we need to introduce some notations.
For $m\in \Z$, we denote by ${\cal A}_m$ the space of continuous functions  
$a:\R_+\rightarrow \C$ such that\\
\noindent i) $a\in C^\infty(\R_+^*)$, $\supp a\subset \{r\leq 2\delta\}$;\\
\noindent ii) for $0\leq r< \delta$, $a$ has an absolutely convergent expansion of the form
$$a(r)=\sum\limits_{{n\geq K(m)\atop n-m-1\in 2\Z}} \sum\limits_{l=0}^n\alpha_{n,l}(\ln r)^lr^{2\nu n},$$
where  $K(m)=m+1$ if $m\geq 0$, and $K(m)=|m|-1$ if $m\leq -1$.
For $k\geq 1$ we define ${\cal B}_k$ as the space of continuous functions  
$b:\R_+\rightarrow \C$ such that\\
\noindent i) $b\in C^\infty(\R_+^*)$;\\
\noindent ii) for $0\leq r< \delta$, $b$  has an absolutely convergent expansion of the form
$$b(r)=\sum\limits_{n=0}^\infty \sum\limits_{l=0}^{2n} \beta_{n,l}r^{4\nu n}(\ln r)^l,$$
\noindent iii) for $r\geq 2\delta$, $b$ is a polynome of degree $k-1$.\\
\noindent  Finally, we set ${\cal B}_k^0=\{b\in {\cal B}_k, b(0)=0\}$.

Clearly, for any $m$, $k$, one has $r\partial_r {\mathcal A}_m\subset {\mathcal A}_m$,
$r\partial_r {\mathcal B}_k\subset {\mathcal B}_k$,
${\mathcal B}_k{\mathcal A}_m\subset {\mathcal A}_m$.
 Note also that
\begin{equation*}\begin{split}
&f_0\in r^{2i\alpha_0}{\cal A}_{0},\quad\varphi\in {\cal B}_1^0,\quad
g_{1,0,0,0}\in r^{2i\alpha_0-2}{\cal A}_{0},\\
& g_{1,2j+1,-1,s}\in r^{-2i\alpha_0-2\nu(2j+1)-2}{\cal B}_{1},\quad 0\leq s\leq 2j+1.
\end{split}
\end{equation*}

Furthermore, one checks easily that if for all $(k,q, m,s)\in \Omega$,
$g_{k, q, m, s}\in r^{2i\alpha_0(1+2m)-2\nu q -2k} {\cal A}_{m}$ if $m\neq -1$ and 
$g_{k, q, -1, s}\in  r^{-2i\alpha_0-2\nu q -2k} {\cal B}_{k}$, then
\be\label{structure}\begin{split}
&\Psi^{lin, i}_{k,q,m,s}, \Psi^{nl,j}_{k,q,m,s}, \tilde \Psi^{nl, 2}_{k,q,m,s} \in 
 r^{2i\alpha_0(1+2m)-2\nu q -2(k-1)} {\cal A}_{m}, \quad m\neq  -1,\\
&\Psi^{lin, 2}_{k,q,-1,s}, \Psi^{nl,j}_{k,q,-1,s}, \tilde \Psi^{nl, 2}_{k,q,-1,s} \in 
 r^{-2i\alpha_0-2\nu q -2(k-1)} {\cal B}_{k-2}, 
\end{split}
\ee
$i=1,2,\,\,j=0,1,2$.

Consider  \eqref{R2}.
Using \eqref{remote:01},  \eqref{remote:10},  \eqref{remote:100},  \eqref{N2}, \eqref{phi},
one can rewrite it as
\be\label{R21}\left\{
\begin{array}{ll}
&\frac14m(m+1)r^2g_{k,q,m,s}=B_{k,q,m,s},\quad m\neq 0, -1,\\
&r\partial_rg_{k,q, m,s}+
\left (2\nu q+k+1+2i\a_0-\frac{r(\bar f_0\partial_rf_0+ f_0\partial_r\bar f_0)}{1+|f_0|^2}\right)g_{k,q, -1,s}=C_{k,q, -1,s},\\
&(2\nu q+k)g_{k,q,0,s}-(s+1)g_{k,q,0,s+1}=C_{k,q, 0,s}+D_{k,q,s},
\end{array}\right.\ee
 where $B_{k,q,m,s}$, $C_{k,q,m,s}$ depend on $g_{k^\prime, q^\prime, m^\prime, s^\prime}$, 
$k^\prime\leq k-1$ only:
\begin{equation*}\begin{split}
&B_{k,q,m,s}=B_{k,q,m,s}(r; g_{k^\prime, q^\prime, m^\prime, s^\prime},\, k^\prime\leq k-1),
\quad m\neq 0,-1,\\
&C_{k,q,m,s}=C_{k,q,m,s}(r; g_{k^\prime, q^\prime, m^\prime, s^\prime},\, k^\prime\leq k-1),
\quad m= 0,-1,\end{split}\end{equation*}
and have the following form
\begin{equation}\label{BC1}\begin{split}
&B_{k,q,m,s}=-\Psi^{lin,1}_{k,q,m,s}-\Psi^{lin,2}_{k,q,m,s}-\Psi^{nl}_{k,q,m,s},\quad m\neq 0,-1\\
&C_{k,q,m,s}=-i\Psi^{lin,2}_{k+1,q,m,s}-i\tilde\Psi^{nl}_{k+1,q,m,s},
\quad m= 0,-1.\end{split}\end{equation}
Finally $D_{k,q,s}$ depend only on $g_{k,q,1,s}$ and is given by
\begin{equation}\label{D}
D_{k,q,s}=-i\Psi^{nl,2,0}_{k+1,q,0,s}=-i\frac{r^2\bar f_0}{1+|f_0|^2}\sum\limits_{{q_1+q_2=q\atop s_1+s_2=s}}g_{1,q_1, -1,s_1}g_{k,q_2, 1,s_2}.
\ee
Note that $D_{2,q,s}=0$.
\begin{remark}\label{rem}
It is not difficult to check that if 
\begin{equation*}\begin{split}
&g_{k,q,m,s}=0,\quad \forall q>(2N+1)(2k-2), \,\, m\neq 0,1,\\
&g_{k,q,m,s}=0,\quad \forall q>(2N+1)(2k-1), \,\, m= 0,1,
\end{split}
\end{equation*}
then 
\begin{equation*}\begin{split}
&B_{k,q,m, s}=0,\quad  \forall q>(2N+1)(2k-2), \,\, m\neq 0,1,\\
&C_{k,q,m, s}=0,\quad  \forall q>(2N+1)(2k-1), \,\, m= 0,1,\\
&D_{k,q,s}=0,\quad  \forall q>(2N+1)(2k-1).
\end{split}
\end{equation*}
\end{remark}
We are now in position to prove the following result.
\begin{lemma}\label{lr}
There exists a unique solution $(g_{k,q,m,s})_{{(k,q,m,s)\in \Omega\atop k\geq 2}}$ of
\eqref{R21} verifying
\be\label{lr1}\begin{split}
&g_{k,q,m,s}\in r^{2i\alpha_0(2m+1)-2\nu q-2k}{\mathcal A}_m,\quad m\neq -1,\\
&g_{k,q,-1,s}\in r^{-2i\alpha_0-2\nu q-2k}{\mathcal B}_k.
\end{split}
\ee
In addition, one has  
\begin{equation}\label{lr2}\begin{split}
&g_{k,q,m,s}=0,\quad \forall q>(2N+1)(2k-2), \,\, m\neq 0,1,\\
&g_{k,q,m,s}=0,\quad \forall q>(2N+1)(2k-1), \,\, m= 0,1,
\end{split}
\ee
\end{lemma}
{\it Proof}. For $k=2$ \eqref{R21}, \eqref{BC1},  \eqref{N2} give
\begin{align}
&\frac12r^2g_{2,2j,-2,s}=B_{2,2j,-2,s},\quad 0\leq s\leq 2j,\, \, 1\leq j,\label{k=2;1}\\
\begin{split}
&r\partial_rg_{2,2j+1, -1,s}+
\left (2\nu (2j+1)+3+2i\a_0-\frac{r(\bar f_0\partial_rf_0+ f_0\partial_r\bar f_0)}{1+|f_0|^2}\right)
g_{2,2j+1, -1,s}\\
&=C_{2,2j+1, -1,s},\quad 0\leq s\leq 2j+1,\, \,0\leq j, \end{split}\label{k=2;2}\\
&(4\nu j+2)g_{2,2j,0,s}-(s+1)g_{2,2j,0,s+1}=C_{2,2j, 0,s},\quad 0\leq s\leq 2j,\, \, 0\leq j,\label{k=2;3}
\end{align}
Recall that $B_{2,q, m, s}$, $C_{2,q, m,s}$ depend only on $g_{1, q^\prime, m^\prime, s^\prime}$
and therefore, are known  by now. By   \eqref{structure}, \eqref{BC1} and remark \ref{rem}
they verify
\begin{equation*}\begin{split}
&B_{2,q, -2, s}\in 
r^{-6i\alpha_0-2\nu q -2} {\cal A}_{-2}, \quad m\neq  0,-1\\
&C_{2,q, 0, s}\in
r^{2i\alpha_0-2\nu q -4} {\cal A}_{0},\quad
C_{2,q, -1, s}\in  r^{-2i\alpha_0-2\nu q -4}{\cal B}_{1},\\
&B_{2,q, -2, s}=0,\quad q>2(2N+1),\\
&C_{2,q, m, s}=0, \quad q>3(2N+1),\,\, m=0,1.
\end{split}
\end{equation*}
Therefore, 
 we get from \eqref{k=2;1}, \eqref{k=2;2},
\be\label{k=2;4}\begin{split}
&g_{2,2j,-2,s}=\frac{2}{r^2}B_{2,2j,-2,s}\in r^{-6i\alpha_0-4\nu j-4} {\cal A}_{-2},\quad 0\leq s\leq 2j,\, \, 1\leq j,\\
&g_{2, 2j, 0, 2j}=\frac{1}{4j\nu+2}C_{2,2j, 0,2j}\in r^{2i\alpha_0-4\nu j-4} {\cal A}_{0},\,\,0\leq j,\\
&g_{2, 2j, 0, s}=\frac{1}{4j\nu+2}C_{2,2j, 0,s}+\frac{s+1}{4j\nu+2}g_{2, 2j, 0, s+1} \in r^{2i\alpha_0-4\nu j-4} {\cal A}_{0},\,\,  0\leq s\leq 2j,\\
&g_{2,2j,-2,s}=0,\quad j>2N+1,\\
&g_{2, 2j, 0, s}=0, \quad j\geq 3N+2,\\
\end{split}
\end{equation}
Consider  \eqref{k=2;3}. Write
$$ g_{2, 2j+1, -1, s}=r^{-2i\alpha_0-3-2\nu(2j+1)}(1+|f_0|^2)\hat  g_{2, 2j+1, -1, s}.$$
Then $\hat  g_{2, 2j+1, -1, s}$ solves
\be\label{k=2;5}
\partial_r\hat g_{2,2j+1, -1,s}=r^{-2}\hat C_{2,2j+1, -1,s},
\ee
where
$$\hat C_{2,2j+1, -1,s}=r^{2i\alpha_0+4+2\nu(2j+1)}(1+|f_0|^2)^{-1} C_{2,2j+1, -1,s}.$$
Since 
$C_{2,2j+1, -1, s}\in  r^{-2i\alpha_0-2\nu (2j+1) -4}{\cal B}_{1}$, we have:\\
 (i) for $0\leq r< \delta$, $\hat C_{2,2j+1, -1,s}$ admits an absolutely convergent  expansion
of the form
$$\hat C_{2,2j+1, -1,s}=\sum\limits_{n=0}^\infty \sum\limits_{l=0}^{2n} \beta_{n,l}r^{4\nu n}(\ln r)^l,$$
 (ii) for $r\geq 2\delta$, $\hat C_{2,2j+1, -1,s} $ is a constant.\\ 
Clearly, there exists a unique solution $\hat g_{2,2j+1, -1,s}$ of \eqref{k=2;5} such that
$\hat g_{2,2j+1, -1,s}\in r^{-1}{\mathcal B}_2$. It is given by
$$\hat g_{2,2j+1, -1,s}(r)=\int_0^r d\rho 
\rho^{-2}(\hat C_{2,2j+1, -1,s}(\rho)-\beta_{0,0})-\beta_{0,0}r^{-1},
\,\, 0\leq s\leq 2j+1,\,\, 0\leq j.$$
Finally, 
since $C_{2,q, -1, s}=0$ for $ q>3(2N+1)$, one has 
$$ g_{2,2j+1, -1,s}=0,\quad j>3N+1.$$

We next proceed by induction.
Suppose we have solved 
\eqref{R21} with $k=2,\dots, l-1$, $l\geq 3$, and have found
$(g_{k,q,m,s})_{{(k,q,m,s)\in \Omega\atop 2\leq k\leq l-1}}$ verifying \eqref{lr1} and \eqref{lr2}.
Consider $k=l$. From the first line in \eqref{R21} we have:
$$\frac14m(m+1)r^2g_{l,q,m,s}=B_{l,q,m,s},\quad m\neq 0, -1,$$
where $B_{l,q,m,s}$ are known by now and, by \eqref{structure}, \eqref{BC1} and remark \ref{rem},
satisfy
\begin{equation*}\begin{split}
&B_{l,q, m, s}\in 
r^{2i\alpha_0(2m+1)-2\nu q -2(l-1)} {\cal A}_{m}, \\
&B_{l,q, m, s}=0,\quad q>2(2N+1)(2l-2).
\end{split}
\end{equation*}
As a consequence, one obtains for $m\neq 0, -1$:
\begin{equation}\begin{split}\label{***}
&g_{l,q,m,s}=\frac{4}{m(m+1)r^2}B_{l,q, m, s}\in r^{2i\alpha_0(2m+1)-2\nu q -2l} {\cal A}_{m},\\
&g_{l,q,m,s}=0,\quad  q>2(2N+1)(2l-2).
\end{split}
\end{equation}
We next consider the equations  for $g_{l,2j,0,s}$:
\begin{equation}\label{k=l, m=0}
(4\nu j+l)g_{l,2j,0,s}-(s+1)g_{l,2j,0,s+1}=C_{l,2j, 0,s}+D_{l,2j,s},
\quad 0\leq s\leq 2j,\,\,0\leq j.
\ee
The right hand side $C_{l,2j, 0,s}+D_{l,2j,s}$ depends only on $g_{l,q_1, 1, s_1}$ and 
$g_{k, q_2, m_2, s_2}$, $k\leq l-1$, and by  \eqref{structure}, \eqref{BC1},
\eqref{***}  and remark \ref{rem},
satisfies
\begin{equation*}\begin{split}
&C_{l,2j, 0,s}+D_{l,2j,s}\in r^{2i\alpha_0-4\nu j -2l} {\cal A}_{0}, \\
&C_{l,2j, 0,s}+D_{l,2j,s}=0,\quad j>(2N+1)(2l-1).
\end{split}
\end{equation*}
Therefore, the solution of  \eqref{k=l, m=0} verifies
\begin{equation*}\begin{split}
&g_{l,2j,0,s}\in r^{2i\alpha_0-4\nu j -2l} {\cal A}_{0},\quad 0\leq s\leq 2j, \,\, 0\leq j,\\
&g_{l,2j,0,s}=0,\quad j>(2N+1)(2l-1).
\end{split}
\end{equation*}
Finally for $g_{l,2j+1,-1,s}$, $0\leq s\leq 2j+1$, $0\leq j$ we have
\begin{equation}\label{k=l, m=-1}\begin{split}
r\partial_rg_{l,2j+1, m,s}&+
\left (2\nu (2j+1)+l+1+2i\a_0-\frac{r(\bar f_0\partial_rf_0+ f_0\partial_r\bar f_0)}{1+|f_0|^2}\right)g_{l,2j+1, -1,s}\\
&=C_{l,2j+1, -1,s},
\end{split}\ee
with $C_{l,2j+1, -1,s}\in r^{-2i\alpha_0-2\nu (2j+1)-2l}{\mathcal B}_{l-1}$ such that
\begin{equation}\label{****}
C_{l,2j+1, -1,s}=0,\quad 2j+1>(2N+1)(2l-1).\ee
Equation \eqref{k=l, m=-1} has a unique solution $g_{l,2j+1,-1,s}$ verifying 
$g_{l,2j+1,-1,s}\in  r^{-2i\alpha_0-2\nu (2j+1)-2l}{\mathcal B}_{l}$, which is given by
\begin{equation*}\begin{split}
g_{l,2j+1,-1,s}=&r^{-2i\alpha_0-2\nu (2j+1)-l-1}(1+|f_0|^2)\hat g_{l,2j+1,-1,s},\\
\hat g_{l,2j+1,-1,s}=&
\int_0^r d\rho\rho^{-l}\big(\hat C_{l,2j+1, -1,s}-
\sum_{0\leq n\leq \frac{l-1}{4\nu}}\sum_{p=0}^{2n}\beta_{n,p}\rho^{4\nu n}(\ln \rho)^p \big)\\
&-
\int_r^\infty d\rho\rho^{-l}\sum_{0\leq n\leq \frac{l-1}{4\nu}}\sum_{p=0}^{2n}\beta_{n,p}\rho^{4\nu n}(\ln \rho)^p ,
\end{split}
\end{equation*}
where
\begin{equation*}\begin{split}
&\hat C_{l,2j+1, -1,s}= r^{2i\alpha_0+2\nu (2j+1)+2l}(1+|f_0|^2)^{-1}C_{l,2j+1, -1,s},\\
&\hat C_{l,2j+1, -1,s}=\sum_{n=0}^\infty\sum_{p=0}^{2n}\beta_{n,p}r^n(\ln r)^p,\quad r<\delta.
\end{split}
\end{equation*}
By \eqref{****},
$$g_{l,2j+1,-1,s}=0, \quad 2j+1>(2N+1)(2l-1).\quad\quad\quad \square$$

Let us define
\begin{equation*}\begin{split}
&w_{rem}^{(N)}(r,t)=f_0(r)+\sum_{(k,q,m,s)\in \Omega,\, k\leq N}t^{k+2\nu q}e^{-im\Phi}(\ln r-\ln t)^sg_{k,q, m,s}(r),\\
&A_{rem}^{(N)}=-i\partial_tw_{rem}^{(N)}-\Delta w_{rem}^{(N)}+r^{-2}w_{rem}^{(N)} +
G(w_{rem}^{(N)},\bar w_{rem}^{(N)},\partial_r w_{rem}^{(N)})\\
&W_{rem}^{(N)}(y,t)=e^{-i\alpha(t)}w_{rem}^{(N)}(rt^{-1/2},t).
\end{split}
\end{equation*}
As a direct consequence of the previous analysis we get:
\begin{lemma}\label{rem.l}
There exists $T(N,\delta)>0$ such that for
 $0<t\leq T(N,\delta)$ the following holds.\\
(i) For any $0\leq l,k\leq 4$, $i=0,1$ and  $ \frac{1}{10}t^{-\varepsilon_2}\leq  y \leq {10}t^{-\varepsilon_2}$,
one has
\be\label{ss.1***}
|y^{-l}\partial_y^{k}\partial_t^i(W_{ss}^{(N)}-W_{rem}^{(N)})|
\leq t^{\nu(1-2\varepsilon_2)N}+t^{\varepsilon_2N},
\ee
provided $N$ is sufficiently large (depending on $\varepsilon_2$).\\
(ii) The profile $w_{rem}^{(N)}(r,t)$ verifies
\begin{align}
&\|r^{-l}\partial_r^k (w_{rem}^{(N)}(t)-f_0)\|_{L^2(rdr,  r\geq \frac{1}{10}t^{1/2-\varepsilon_2})}
\leq Ct^{\eta},\quad 0\leq k+l\leq 3,\label{rem.1.11}\\
&\|r\partial_rw_{rem}^{(N)}(t)\|_{L^\infty(r\geq \frac{1}{10}t^{1/2-\varepsilon_2})}
\leq C\delta^{2\nu},\label{rem.1.111}\\
&\|r^{-l}\partial_r^kw_{rem}^{(N)}(t)\|_{L^\infty(r\geq \frac{1}{10}t^{1/2-\varepsilon_2})}
\leq C(\delta^{2\nu-k-l}+ t^{\nu-(k+l)/2+\eta}),\quad 0\leq k+l\leq 4,\label{rem.1.2}\\
&\|r^{-l-1}\partial_r^kw_{rem}^{(N)}(t)\|_{L^\infty(r\geq \frac{1}{10}t^{1/2-\varepsilon_2})}
\leq C(\delta^{2\nu-6}+ t^{\nu-3+\eta}),\quad k+l=5\label{rem.1.12}
\end{align}
(iii) The error $A_{rem}^{(N)}(r,t)$ admits the estimate
\be\label{rem.1.21}
\|r^{-l}\partial_r^k\partial_t^iA_{rem}^{(N)}(t)\|_{L^2(rd r, r\geq \frac{1}{10}t^{1/2-\varepsilon_2})}
\leq t^{\varepsilon_2N},\quad 0\leq l+k\leq 3,\,\, i=0,1,
\ee
provided $N$ is sufficiently large.
\end{lemma}
\subsubsection{Proof of proposition \ref{p1}}
We are now in position to finish the proof of  proposition \ref{p1}.
Fix $\varepsilon_2$ verifying $0<\varepsilon_2<\frac12$.
For $N\geq 2$, define
\begin{equation*}\begin{split}
\hat W_{ex}^{(N)}(\rho,t)=&\theta(t^{\nu-\varepsilon_1}\rho)W_{in}^{(N)}(t^{\nu}\rho,t)+
(1-\theta(t^{\nu-\varepsilon_1}\rho))\theta(t^{\nu+\varepsilon_2}\rho)W_{ss}^{(N)}(t^{\nu}\rho,t)\\
&+(1-\theta(t^{\nu+\varepsilon_2}\rho))e^{-i\alpha(t)}w_{rem}^{(N)}(t^{\nu+1/2}\rho,t),\\
V_{ex}^{(N)}(\rho,t)=&\big(\frac{2\re  \hat W_{ex}^{(N)}}{1+|\hat W_{ex}^{(N)}|^2}, \frac{2\im  
\hat W_{ex}^{(N)}}{1+| \hat W_{ex}^{(N)}|^2},
\frac{1-| \hat W_{ex}^{(N)}|^2}{1+|\hat  W_{ex}^{(N)}|^2}).
\end{split}
\end{equation*}
Clearly, $V^{(N)}_{ex}(\rho, t)$ is well defined for $\rho$ is sufficiently large,
and for $ \rho<t^{-\nu+\varepsilon_1}$ $V^{(N)}_{ex}(\rho, t)$ coincides with
$V^{(N)}_{in}(\rho,t)$.
Therefore, setting
\begin{equation*}\begin{split}
&V^{(N)}(\rho,t)=\left\{\begin{array}{cr}& V^{(N)}_{in}(\rho,t), \quad \rho\leq \frac12t^{-\nu+\varepsilon_1},
\\&V^{(N)}_{ex}(\rho,t)\quad \rho\geq \frac12t^{-\nu+\varepsilon_1}.\end{array}\right.\\
&u^{(N)}(x,t)=e^{(\alpha(t)+\theta)R}V^{(N)}(\lambda(t)|x|,t),
\end{split}
\end{equation*}
we get a  $C^\infty$ 1- equivariant profile $u^{(N)}: \R^2\times \R_+^*\rightarrow S^2$
that, by lemmas \ref{in} (i), \ref{ss} (ii), \ref{rem.l} (ii),  for any $N\geq 2$ verifies
part (i) of proposition \ref{p1}, $\zeta^*_N$ being given by
$$\zeta^*_N(x)=e^{\theta R}\hat \zeta^*_N(|x|),\quad \hat\zeta^*_N=\big(\frac{2\re  f_0}{1+|f_0|^2}, \frac{2\im f_0}{1+|f_0|^2},
\frac{1-|f_0|^2}{1+|f_0|^2}).$$
By lemmas  \ref{in} (ii), \ref{ss} (i), (iii) and \ref{rem.l} (i), (iii), for $N$ sufficiently large
the  error $r^{(N)}=-u^{(N)}_t+u^{(N)}\times \Delta u^{(N)}$
satisfies
$$\|{ r}^{(N)}(t)\|_{H^{3}}+ \|\partial_t{ r}^{(N)}(t)\|_{H^{1}}+
\|<x>r^{(N)}(t)\|_{L^{2}}
\leq t^{\eta N},\quad t\leq T(N,\delta),$$
with some $\eta=\eta(\nu, \varepsilon_2)>0$.
Re-denoting  $N=\frac{N}{\eta}$ we obtain a family of approximate solutions $u^{(N)}(t)$ verifying  proposition \ref{p1}.
\section{Proof of the theorem} 
\subsection{Main proposition}
The proof of theorem 1.1 will be achieved  by compactness arguments that rely on the following 
result.
Let $u^{(N)}, T=T(N,\delta)$ be as in proposition 2.1. 
 Consider the Cauchy problem
\be\label{3.1}\begin{split}
&u_t=u\times \triangle u, \quad t\geq t_1,\\
&u|_{t=t_1}=u^{(N)}(t_1),\end{split}\ee
with $0<t_1<T$. 

One has
\begin{proposition}\label{mp}
For $N$ sufficiently large there exists $0<t_0<T$ such that for any $t_1\in(0, t_0)$ the solution $u(t)$ of
\eqref{3.1} verifies:\\
\noindent  (i) $u-u^{(N)}$ is in $C([t_1,t_0], H^3)$ and one has 
\be\label{3.1*}
\|u(t)-u^{(N)}(t)\|_{H^3}\leq t^{N/2}, \quad\forall t_1\leq t\leq t_0.\ee
\noindent  (ii ) Furthermore, $<x>(u(t)-u^{(N)}(t))\in L^2$ and
\be\label{3.4}
\|<x>(u(t)-u^{(N)}(t))\|_{L^2}\leq t^{N/2},\quad\forall t_1\leq t\leq t_0.
\ee
\end{proposition}
{\it Proof}. The proof is by bootstrap argument.
Write
$$
u^{(N)}(x,t)=e^{\alpha(t) R}U^{(N)}(\lambda( t)x,t), \quad r^{(N)}(x,t)=\lambda^2(t)e^{\alpha(t) R}R^{(N)}(\lambda( t)x,t)
$$
$$
u(x,t)=e^{\alpha(t) R}U(\lambda (t)x,t),
\quad U(y,t)=U^{(N)}(y,t)+S(y,t),$$
$$U^{(N)}(y,t)=\phi(y)+\chi^{(N)}(y,t).$$
Then $S(t)$ solves
\be\label{3.10}
t^{1+2\nu}S_t+\alpha_0t^{2\nu}RS-(\nu+{1\over 2})y\cdot\nabla S=
S\times \Delta U^{(N)}+ U^{(N)}\times \Delta S+S\times \Delta S +R^{(N)}(t).
\ee
Assume that
\be\label{boot}
\|S\|_{L^\infty(\R^2)}\leq \delta_1,\ee
with $\delta_1$ sufficiently small.
Note that since $S$ is 1-equivariant and
\be\label{E0}
(\phi, S)+(\chi^{(N)}, S)+|S|^2=0
\ee
where  
$\|\chi^{(N)}\|_{L^\infty(\R^2)}\leq C\delta^{2\nu}$ (see \eqref{p1.1.1}),
the bootstrap assumption \eqref{boot} implies
\be\label{boot.1}
\|S\|_{L^\infty(\R^2)}\leq C\|\nabla S\|_{L^\infty(\R^2)}.
\ee
\subsubsection{Energy control}
We will first derive a bootstrap control of the energy norm:
$$J_1(t)=\int_{\R^2} dy (|\nabla S|^2+\kappa(\rho)|S|^2),\quad \rho=|y|.$$
It follows from \eqref{3.10} that
\begin{equation*}\begin{split}
&t^{1+2\nu}\frac{d}{dt}J_1(t)={\mathcal E}_1+{\mathcal E}_2+{\mathcal E}_3+{\mathcal E}_4,\\
&{\mathcal E}_1=-2\int dy (S\times \Delta \chi^{(N)},\Delta S),\\
&{\mathcal E}_2=2\int dy \kappa(\chi^{(N)}\times \Delta S,S),\\
&{\mathcal E}_3=-(\frac12 +\nu)t^{2\nu}\int dy (2\kappa+\rho\kappa^\prime)( S,S),\\
&{\mathcal E}_4=2\int dy\big[(\nabla R^{(N)},\nabla S) +\kappa ( R^{(N)},S)\big].
\end{split}
\end{equation*}
From proposition \ref{p1} we have
\begin{equation*}\begin{split}
|{\mathcal E}_j|&\leq Ct^{2\nu}\|S\|^2_{H^1},\quad j=1,\dots, 3,\\
|{\mathcal E}_4|&\leq Ct^{N+\nu +1/2}\|\nabla S\|_{L^2}.
\end{split}
\end{equation*}

Combining these inequalities we obtain
\be\label{E8}
\big|\frac{d}{dt}J_1(t)\big|\leq Ct^{-1}\|S\|^2_{H^1}+
Ct^{2N-1/2-3\nu}.
\ee
\subsubsection{Control of the $L^2$ norm}
Consider $J_0(t)=\int_{\R^2} dy |S|^2$. We have
\begin{equation*}\begin{split}
&t^{1+2\nu}\frac{d}{dt}J_0(t)={\mathcal E}_{5}+{\mathcal E}_{6}+{\mathcal E}_{7},\\
&{\mathcal E}_{5}=2\int dy (U^{(N)}\times \Delta S,S),\\
&{\mathcal E}_{6}=-2(1+2\nu)t^{2\nu}J_0(t),\\
&{\mathcal E}_{7}=2\int dy (R^{(N)},S).
\end{split}
\end{equation*}
Consider ${\mathcal E}_{5}$. Decomposing $U^{(N)}$ and $S$ in the basis
$f_1$, $f_2$, $Q$:
$$U^{(N)}(y,t)=e^{\theta R}((1+z_3^{(N)}(\rho,t))Q(\rho)+z_1^{(N)}(\rho,t)f_1(\rho)+z_2^{(N)}(\rho,t)f_2(\rho)),$$
$$S(y,t)=e^{\theta R}(\zeta_1(\rho,t)f_1(\rho)+\zeta_2(\rho,t)f_2(\rho)+\zeta_3(\rho,t)Q(\rho)),$$
one can rewrite ${\mathcal E}_{5}$ as follows.
\begin{equation*}\begin{split}
&{\mathcal E}_{5}={\mathcal E}_{8}+{\mathcal E}_{9}+{\mathcal E}_{10},\\
&{\mathcal E}_{8}=-4\int _{\R_+}d\rho\rho\frac{ h_1}{\rho}\zeta_2\partial_\rho\zeta_3,\\
&{\mathcal E}_{9}=-2\int _{\R_+}d\rho\rho(\partial_\rho z^{(N)}\times\partial_\rho  \zeta, \zeta),
\quad z ^{(N)}=(z_1^{(N)},z^{(N)}_2, z_3^{(N)}),\, \zeta=(\zeta_1, \zeta_2,\zeta_3),\\
&{\mathcal E}_{10}=2\int _{\R_+}d\rho \rho(z^{(N)}\times l,\zeta),
\end{split}
\end{equation*}
where
\begin{equation*}
l=(-\frac{1}{\rho^2}\zeta_1-\frac{2h_1}{\rho}\partial_\rho\zeta_3,-\frac{1}{\rho^2}\zeta_2,
\kappa(\rho)\zeta_3+\frac{2h_1}{\rho}\partial_\rho\zeta_1-\frac{2h_1h_3}{\rho^2}\partial_\rho\zeta_1).
\end{equation*}
Clearly,
$$|l|\leq C\rho^{-2}(|\zeta|+|\partial_\rho \zeta|).$$
Therefore,
\be\label{M1}
|{\mathcal E}_{10}|\leq Ct^{2\nu}\|S\|^2_{H^1}.
\ee
Consider ${\mathcal E}_{8}$. It follows from 
\be\label{M0}
2(\zeta,\vec k+z^{(N)})+|\zeta|^2=0,
\ee that
$$|\partial_\rho\zeta_3|\leq C(|\partial_\rho z^{(N)}||\zeta| +|z^{(N)}||\partial_\rho\zeta|+
|\partial_\rho\zeta||\zeta|).$$
As a consequence,
\be\label{M2}
|{\mathcal E}_{8}|\leq C\big[t^{2\nu}\|S\|^2_{H^1}+
\|\nabla S\|^3_{L^2}\big].
\ee
Consider ${\mathcal E}_{9}$. 
Denote $e_0=\vec k+ z^{(N)}$ 
and write $\zeta=\zeta^{\perp}+\mu e_0$, $\mu=(\zeta,e_0)$.
It follows from \eqref{M0} that 
\begin{equation*}\begin{split}
&|\mu|\leq C|\zeta|^2,\\
&|\mu_\rho|\leq C|\zeta||\partial_\rho \zeta|.
\end{split}
\end{equation*}
Therefore,
 ${\mathcal E}_{9}$ can be written as
\be\label{M3}
{\mathcal E}_{9}=-2\int _{\R_+}d\rho\rho(\partial_\rho\zeta^{\perp}\times \zeta^{\perp}, \partial_\rho e_0) +O(\|S\|_{H^{1}}^2\|\nabla S\|_{L^2}).
\ee
Let $e_1,\, e_2$ be a smooth orthonormal basis of the tangent space $T_{e_0}S^2$
that verifies $e_2=e_0\times e_1$. Then
the expression
$(\partial_\rho\zeta^{\perp}\times \zeta^{\perp}, \partial_\rho e_0)$ can be written as follows:
$$
(\partial_\rho\zeta^{\perp}\times \zeta^{\perp}, \partial_\rho e_0)=
(\zeta^{\perp}, \partial_\rho e_0)\left[(\zeta^{\perp}, e_2)(\partial_\rho e_0,e_1)-
(\zeta^{\perp}, e_1)(\partial_\rho e_0,e_2)\right],
$$
which leads to the estimate
\be\label{M4}
\left|\int _{\R_+}d\rho\rho(\partial_\rho\zeta^{\perp}\times \zeta^{\perp}, \partial_\rho e_0)\right|
\leq C\|\partial_\rho z^{(N)}\|_{L^\infty}^2J_0(t)\leq Ct^{2\nu}J_0(t).
\ee
Combining \eqref{M1}, \eqref{M2}, \eqref{M3}, \eqref{M4} we obtain
\be\label{M5}
\big|\frac{d}{dt}J_0(t)\big|\leq C\left[t^{-1}\| S\|_{H^1}^2+
t^{-1-2\nu}\|S\|_{H^{1}}^2\|\nabla S\|_{L^2} +t^{2N-1/2-3\nu}\right].
\ee

\subsubsection{Control of the weighted  $L^2$ norm}
Using \eqref{3.10} to compute the derivative $\frac{d}{dt}\|yS(t)\|_{L^2}^2$,
we obtain
\begin{equation*}\begin{split}
t^{1+2\nu}\frac{d}{dt}\||y|S(t)\|_{L^2}^2=&-4\int dy y_i(U^{(N)}\times \partial_iS, S)\\
&-2\int dy |y|^2(\partial_iU^{(N)}\times \partial_iS, S)\\
&-2(1+2\nu)t^{2\nu}\||y|S(t)\|_{L^2}^2
+2\int dy|y|^2(R^{(N)},S).
\end{split}
\end{equation*}
Here and below $\partial_j$ stands for $\partial_{y_j}$, the summation over  the repeated indexes
being  assumed.

As a consequence, we get
\be\label{MW11}
\left|\frac{d}{dt}\||y|S(t)\|_{L^2}^2\right|\leq \frac{C}{t}\left[
\||y|S(t)\|_{L^2}^2+t^{-4\nu}\|S\|_{H^{1}}^2+t^{2N-2\nu}\right].
\ee
\subsubsection{Control of the higher regularity}
In addition to \eqref{boot},  assume that
\be\label{boot1}
\|S(t)\|_{H^3}+\| |y|S(t)\|_{L^2}\leq t^{2N/5}.
\ee
We will control $\dot H^3$ norm of the solution by means of 
$\|\nabla S_t\|_{L^2}$. More precisely, consider the functional
$$J_3(t)=t^{2+4\nu}\int dx |\nabla s_t(x,t)|^2+ t^{1+2\nu}\int dx \kappa(t^{-1/2-\nu}x)|s_t(x,t)|^2,$$
where $s(x,t)$ is defined by
$$s(x,t)=e^{\alpha(t)R}S(\lambda(t)x,t).$$
Write $s_t(x,t)=e^{\alpha(t)R}\lambda^2(t)g(\lambda(t)x,t)$. 
In terms of $g$, $J_3$ can be written as
$J_3(t)=\int dy |\nabla g(y,t)|^2+\int dy \kappa(\rho)|g(y,t)|^2$. Let us compute
the derivative $\frac{d}{dt}J_3(t)$. 
Clearly, $g(y,t)$ solves
\be\label{HE.1}\begin{split}
&t^{1+2\nu}g_t+\alpha_0t^{2\nu}Rg-(\nu+{1\over 2})t^{2\nu}(2+y\cdot\nabla )g=\\
&(S+ U^{(N)})\times \Delta g+g\times ( \Delta U^{(N)}+\Delta S)\\
&+(U^{(N)}\times \Delta U^{(N)}-R^{(N)})\times \Delta S\\
&+S\times \Delta(U^{(N)}\times \Delta U^{(N)}-R^{(N)})+
t^{2+4\nu}r^{(N)}_t.
\end{split}
\ee
Therefore, we get
\be\label{HE.2}
\begin{split}
t^{1+2\nu}\frac{d}{dt}J_3(t)=&(2+4\nu) t^{2\nu} \|\nabla g\|_{L^2}^2+(\frac12+\nu)t^{2\nu}\int  (2\kappa-
\rho\kappa^\prime)|g|^2dy\\
&+
E_1+E_2+E_3+E_4+E_5,
\end{split}
\ee
where
$$
{E}_1=-2\int dy (g\times \Delta \chi^{(N)},\Delta g)+
2\int dy \kappa(\chi^{(N)}\times \Delta g,g),$$

\begin{equation*}\begin{split}
E_2=&-2\int dy((U^{(N)}\times\Delta U^{(N)}-R^{(N)})\times \Delta S, \Delta g)\\
&+
2\int dy(\Delta(U^{(N)}\times\Delta U^{(N)}-R^{(N)})\times S, \Delta g)\\
&+2\int dy\kappa((U^{(N)}\times\Delta U^{(N)}-R^{(N)})\times \Delta S,  g)\\
&-2\int dy\kappa (\Delta(U^{(N)}\times\Delta U^{(N)}-R^{(N)})\times S, g),\\
E_3=&-2\int dy(g\times \Delta S,\Delta g),\\
E_4=&2\int dy\kappa(S\times \Delta g, g),\\
E_5=&-2t^{2+4\nu}\int dy (r_t, \Delta g)+2t^{2+4\nu}\int dy \kappa(r_t,  g).
\end{split}
\end{equation*}
The terms $E_j$, $j=1,4,5$ can be estimated as follows.
\begin{equation}\label{HE.41111}\begin{split}
&|E_1|\leq  Ct^{2\nu}\| g\|_{H^1}^2,\\
&|E_4|\leq  C\| g\|_{H^1}^2\|S\|_{H^{3}}\leq Ct^{2\nu}\| g\|_{H^1}^2,\\
&|E_5|\leq  C(t^{2\nu}\|g\|_{H^1}^2+t^{2N+3+4\nu}),
\end{split}
\end{equation}
provided $N$ is sufficiently large and $t\leq t_0$ with some $t_0=t_0(N)>0$.

For $E_2$ we have 
\begin{equation*}\begin{split}
|E_2|&\leq  C
(\|\Delta\chi^{(N)}\|_{W^{2,\infty}}+\|R^{(N)}\|_{H^3})\| g\|_{H^1}\|S\|_{H^{3}}\\
&+
C\|<y>^{-1}\nabla\Delta^2\chi^{(N)}\|_{L^{\infty}}\|\nabla g\|_{L^2}\|<y>S\|_{L^2}.
\end{split}
\end{equation*}
As a consequence,
\be\label{HE.512}
|E_2|\leq  Ct^{2\nu}(\| g\|_{H^1}\|S\|_{H^{3}}+
\|\nabla g\|_{L^2}\|<y>S\|_{L^2}).
\ee
Note that since 
\be\label{HE.410}
g=(U^{(N)}+S)\times\Delta S+S\times\Delta U^{(N)}+R^{(N)},
\ee
the bootstrap assumption \eqref{boot1} implies
\be\label{HE.4100}\begin{split}
&\|g\|_{L^2}\leq C(\|S\|_{H^{2}}+\|R^{(N)}\|_{L^2}),\\
&\|\nabla g\|_{L^2}\leq C(\|S\|_{H^{3}}+\|\nabla R^{(N)}\|_{L^2}).
\end{split}
\ee
Therefore, \eqref{HE.41111}, \eqref{HE.512} can be rewritten as
\begin{equation}
\label{HE.4}\begin{split}
|E_1|+|E_2|+|E_4|+|E_5|\leq &Ct^{2\nu}[\|S\|_{H^{3}}^2+ 
(\|S\|_{H^{3}}+t^{N+1+2\nu})\|<y>S\|_{L^2}]\\
&+Ct^{2N+1+4\nu}.\end{split}
\ee
Consider $E_3$. One has
\begin{equation}\label{HE.41}\begin{split}
g\times \Delta S&=(U^{(N)}+S,\Delta S)\Delta S-|\Delta S|^2(U^{(N)}+S)\\
&+
(S\times\Delta U^{(N)}+R^{(N)})\times \Delta S,\\
\Delta g=&(U^{(N)}+S)\times\Delta^2 S+Y,\\
Y=2&(\partial_jU^{(N)}+\partial_jS)\times \Delta\partial_j S+S\times \Delta^2 U^{(N)}\\
&+2\partial_jS\times \Delta\partial_jU^{(N)}+\Delta R^{(N)}.
\end{split}
\end{equation}

Therefore, one can write $E_3$ as 
$
E_3=E_6+E_7+E_8,
$
where
\begin{equation*}\begin{split}
E_6=&-2\int dy(U^{(N)}+S,\Delta S)(\Delta S,\Delta g),\\
E_7=&2\int dy|\Delta S|^2(U^{(N)}+S, \Delta g)=2\int dy|\Delta S|^2(U^{(N)}+S,Y),\\
E_8=&-2\int dy((S\times\Delta U^{(N)}+R^{(N)})\times \Delta S,\Delta g).
\end{split}
\end{equation*}
For $E_6$ we have:
\begin{equation*}\begin{split}
E_6=&2\int dy [(\Delta U^{(N)}, S)+2(\partial_jU^{(N)},\partial_j S)+(\partial_jS, \partial_jS)](\Delta S,\Delta g)\\
&=-2\int dy [(\Delta U^{(N)}, S)+2(\partial_jU^{(N)},\partial_j S)+(\partial_jS, \partial_jS)](\Delta\partial_k S,\partial_k g)\\
&-2\int dy(\Delta S,\partial_k g)\partial_k[(\Delta U^{(N)}, S)+2(\partial_jU^{(N)},\partial_j S)+(\partial_jS, \partial_jS)].
\end{split}
\end{equation*}
As a consequence, one obtains:
\begin{equation}\label{HE.10}
|E_6|\leq C\|S\|_{H^3}^2\|g\|_{H^1}\leq Ct^{2\nu}\|S\|_{H^3}^2.
\ee

Consider $E_7$.  From \eqref{HE.41} we have
$$\|Y\|_{L^2}\leq C(\|S\|_{H^3}+t^{N}).$$
Therefore, we obtain:
\begin{equation}\label{HE.12}
|E_7|\leq Ct^{2\nu}\|S\|_{H^3}^2.
\ee

Finally,  $E_8$ can be estimated as follows
\begin{equation}\label{HE.11}
|E_8|\leq C\|g\|_{H^1}(\|S\|_{H^3}^2+t^{N}\|S\|_{H^3})\leq Ct^{2\nu}\|S\|_{H^3}^2+Ct^{3N}.
\ee

Combining \eqref{HE.10}, \eqref{HE.11}, \eqref{HE.12} we get
\begin{equation}\label{HE.13}
|E_3|\leq C(t^{2\nu}\|S\|_{H^3}^2+t^{3N}),
\ee
which together with \eqref{HE.4} gives
\be\label{HE.3}
\left|\frac{d}{dt}J_3(t)\right|\leq \frac{C}{t}\left[\|S\|_{H^{3}}^2+ 
(\|S\|_{H^{3}}+t^{N+1+2\nu})\||y|S\|_{L^2})\right]+Ct^{2N+2\nu}.
\ee
\subsubsection{Proof of proposition \ref{mp}}

To prove the proposition it is sufficient to show that \eqref{boot}, \eqref{boot1} implies \eqref{3.1*}, \eqref{3.4}.

Under the bootstrap assumption \eqref{boot1}, \eqref{E8}, \eqref{M5} become
\be\label{estimate1}
\big|\frac{d}{dt}J_1(t)\big|+\big|\frac{d}{dt}J_0(t)\big|\leq Ct^{-1}\|S\|^2_{H^1}+
Ct^{2N-1/2-3\nu},\quad \forall t\leq t_0,
\ee
provided $N$ is sufficiently large, $t_0$ sufficiently small.

Note that for $c_0>0$ sufficiently large one has
$\|S\|^2_{H^1}\leq J_1+c_0J_0$. Therefore, denoting
$J(t)=J_1(t)+c_0J_0(t)$ one can rewrite \eqref{estimate1} as
\be\label{estimate2}
\big|\frac{d}{dt}J(t)\big|\leq Ct^{-1}J(t)+
Ct^{2N-1/2-3\nu}.
\ee
Integrating this inequality with zero initial condition at $t_1$ one gets 
\be\label{estimate3}
J(t)\leq \frac{C}{N}{t^{2N+1/2-3\nu}},\quad \forall t\in [t_1, t_0],
\ee
provided $N$ is sufficiently large.
As a consequence, we obtain
\be\label{estimate4}
\|S\|^2_{H^1}\leq\frac{C}{N}{t^{2N+1/2-3\nu}},\quad  \forall t\in[ t_1, t_0].
\ee

Consider $\| |y|S(t)\|_{L^2}$. From \eqref{MW11},\eqref{estimate4} we have
\be\label{MW}
\left|\frac{d}{dt}\||y|S(t)\|_{L^2}^2\right|\leq \frac{C}{t}\left[
\||y|S(t)\|_{L^2}^2+t^{2N-7\nu}\right].
\ee
Integrating this inequality and assuming that $N$  is sufficiently large, we get
\be\label{estimate5}
\||y|S(t)\|_{L^2}^2\leq \frac{C}{N}t^{2N-7\nu}, \quad \forall t\in[ t_1, t_0],
\ee
which gives in particular,
\be\label{estimate6}
\||x|s(t)\|_{L^2}^2\leq t^{N/2}, \quad \forall t\in[ t_1, t_0].
\ee
We next consider $\|\nabla\Delta s(t)\|_{L^2(\R^2)}$.
It follows from \eqref{HE.410}, \eqref{boot1} that
for any $j=1, 2$ 
\be\label{est5}
\|\partial_ig-(U^{(N)}+S)\times \Delta \partial_j S\|_{L^2}
\leq C(\| S\|_{H^2(\R^2)}+t^{N+1+2\nu}).
\ee
Note also that since $|U^{(N)}+S|=1$, we have
\begin{equation*}\begin{split}
|(U^{(N)}+S)\times\Delta\partial_jS|^2&=|\Delta\partial_jS|^2-(U^{(N)}+S,\Delta\partial_jS)^2,\\
(U^{(N)}+S,\Delta\partial_jS)=&-(\Delta U^{(N)}+\Delta S,\partial_jS)-\Delta(\partial_jU^{(N)},S)\\
&-2(\partial_kU^{(N)}+\partial_kS,
\partial_{jk}^2S),\end{split}
\end{equation*}
which together with  \eqref{boot1} gives  
\be\label{est6}
\|\Delta\partial_jS\|^2_{L^2}-\|(U^{(N)}+S)\times \Delta \partial_j S\|^2_{L^2}
\leq C\|S\|^2_{H^2}.
\ee

Consider the functional $\tilde J_3(t)=J_3(t)+ c_1J_0(t)$.
It follows from \eqref{HE.4100}, \eqref{est5}, \eqref{est6}   that for $c_1>0$ sufficiently large we have
\be\label{est7}
c_2\| S\|^2_{H^3}-Ct^{2N+1+2\nu} \leq \tilde J_3(t)\leq 
C(\|S\|^2_{H^3}+t^{2N+1+2\nu}),\ee
with some $c_2>0$.

From \eqref{HE.3},  \eqref{estimate1}, \eqref{estimate5} one gets
\be\label{estimate20}
\begin{split}
\big|\frac{d}{dt}\tilde J_3(t)\big|&\leq C\left[t^{-1}(\|S\|_{H^3(\R^2)}^2+\||y|S\|_{L^2(\R^2)}^2)+
t^{2N-1/2-3\nu}\right]\\
&\leq Ct^{-1}\tilde J_3(t)+Ct^{2N-7\nu-1}.
\end{split}
\ee
Integrating this inequality between $t_1$ and $t$ and observing that
$\tilde J_3(t_1)=t_1^{2+4\nu}\int dx |\nabla r^{(N)}(x,t_1)|^2+ t_1^{1+2\nu}\int dx \kappa(t^{-1/2+\nu}x)|r^{(N)}(x,t_1)|^2$, and therefore,
$|\tilde J_3(t_1)|\leq Ct_1^{2N+1+2\nu}$,  we obtain
$$\tilde J_3(t)\leq Ct^{2N-7\nu}, \quad \forall t\in [t_1,t_0].$$
Combining this inequality with \eqref{est7}, one gets
$$\| S\|^2_{H^3(\R^2)}\leq Ct^{2N-7\nu},\quad \forall t\in [t_1,t_0],$$
which implies that
$$\|s\|_{H^3(\R^2)}\leq  t^{N/2},\quad \forall t\in [t_1,t_0].$$
This concludes the proof of  proposition \ref{mp}.

\subsection{Proof of the theorem}
The proof of the theorem is now straightforward. Fix $N$ such that proposition \ref{mp} holds. Take a sequence $(t^j)$,
$0<t^j<t_0$, $t^j\rightarrow 0$ as $j\rightarrow \infty$. Let $u_{j}(x,t)$ be the solution of
\be\label{3.1***}\begin{split}
&\partial_tu_j=u_j\times \Delta u_j,\quad t\geq t^j,\\
&u_j|_{t=t^j}=u^{(N)}(t^j),\end{split}\ee
By proposition \ref{mp}, for any $j$, $u_j-u^{(N)}\in C([t^j, t_0], H^3)$ and satisfies
\be\label{3.11}
\|u_j(t)-u^{(N)}(t)\|_{H^3}+\|<x>(u_j(t)-u^{(N)}(t))\|_{L^2}\leq 2t^{N/2},\quad  \forall t\in [t^j, t_0].\ee
This implies in particular, that  the sequence $u_j(t_0)-u^{(N)}(t_0)$ is compact in $H^2$ and therefore after passing to a subsequence we
can assume that $u_j(t_0)-u^{(N)}(t_0)$ converges in $ H^2$ to  some 1-equivariant function $w\in H^3$,
with  $\|w\|_{H^3}\leq \delta^{2\nu}$, $|u^{(N)}(t_0)+w|=1$.

Consider the Cauchy problem
\be\label{3.12}\begin{split}
&u_t=u\times \Delta u, \quad t\leq t_0,\\
&u|_{t=t_0}=u^{(N)}(t_0)+w,\end{split}\ee
By the  local well-posedness, \eqref{3.12} admits a unique solution $u\in C((t^*,t_0], \dot H^1\cap \dot H^3)$
with some $0\leq t^*<t_ 0$.
By $H^1$ continuity of the flow (see \cite{GKT2}), $u_j\rightarrow u$ in $C((t^{*},t_0], \dot H^1)$,
which together with \eqref{3.11}  gives
\begin{equation}\label{4.1}
\|u(t)-u^{(N)}(t)\|_{H^3}\leq 2t^{N/2},\quad  \forall t\in (t^{*},t_0].\ee
This implies that $t^{*}=0$ and combined with proposition \ref{p1} gives the result stated in 
theorem 1.1.


\begin{thebibliography}{99}
\bibitem{AH} Angenent S.; Hulshof J., Singularities at t = ∞ in equivariant harmonic map flow, Contemp.
Math. 367, Geometric evolution equations, 1–15, Amer. Math. Soc., Providence, RI, 2005.
\bibitem{Ber1} Van den Bergh, J.; Hulshof, J.; King, J., Formal asymptotics of bubbling in the harmonic
map heat flow, SIAM J. Appl. Math. vol 63, o5. pp 1682-1717.
\bibitem{BKT} Bejenaru, I.; Ionescu, A.; Kenig, C.; Tataru, D., Global Schrödinger maps, to appear in Annals
of Math.
\bibitem{BIKT} Bejenaru, I.; Ionescu, A.; Kenig, C.; Tataru, D., Equivariant Schrödinger maps
in two spatial dimensions, arXiv:1112.6122v1.
\bibitem{BT}Bejenaru, I.; Tataru, D., Near soliton evolution for equivariant Schrödinger Maps in two
spatial dimensions, arXiv:1009.1608.
\bibitem{CHU}
Chang, N-H.; Shatah, J.; Uhlenbeck, K., Schrödinger maps, Comm. Pure Appl. Math. 53
(2000), no 5, 590-602.
\bibitem{GS} Grillakis, M.; Stefanopoulos, V., Lagrangian formulation, energy estimates and the Sch\"odinger
map prolem, Comm PDE 27 (2002), 1845-1877.
\bibitem{GS1}Grotowski, J.; Shatah, J. Geometric evolution equations in critical dimensions. Calc. Var.
Partial Differential Equations 30 (2007), no. 4, 499-512.
\bibitem{GKT1} Gustafson, S.; Kang, K.; Tsai, T-P.; Schr\"odinger flow near harmonic maps; Comm. Pure Appl.
Math. 60 (2007), no. 4, 463-499.
\bibitem{GKT2} Gustafson, S.; Kang, K.; Tsai, T-P.; Asymptotic stability of harmonic maps under the
Schr\"odinger flow.; Duke Math. J. 145 no. 3 (2008) 537-583.
\bibitem{GK} Gustafson, S.; Koo, E.; Global well-posedness for $2D$ radial Schr\"odinger maps into
the sphere,arXiv:1105.5659.
\bibitem{GNT}
Gustafson, S.; Nakanishi, K.; Tsai, T-P.; Asymptotic stability, concentration and oscillations
in harmonic map heat flow, Landau Lifschitz and Schr\"odinger maps on R2; Comm. Math.
Phys. (2010), 300, no 1, 205-242.
\bibitem{KST1} Krieger, J.; Schlag, W.; Tataru, D.,  Renormalization and blow up for charge one equivariant
critical wave maps. Invent. Math. 171 (2008), no. 3, 543–615.
 \bibitem{McG} McGahagan, H., An approximation scheme for Schr\"odinger maps, Comm. Partial Differential Equations
32 (2007), 375–40.
\bibitem{MRR} Merle F.; Rapha\"e,l P.; Rodnianski, I.,
Blow up dynamics for smooth equivariant solutions to the energy critical Schrödinger map,  arXiv:1106.0912.
\bibitem{NSU}A. Nahmod, A Stefanov, K. Uhlenbeck, On Schr\"odingers maps, CPAM 56(2003), 114-151.
\bibitem{RR} Raphaël, P.; Rodnianksi, I., Stable blow up dynamics for the critical corotational wave maps
and equivariant Yang Mills problems, to appear in Prep. Math. IHES.
\bibitem{SSB}Sulem, P.L; Sulem, C.; Bardos, C., On the continuous limit for a system of continuous spins,
Comm. Math. Phys 107 (1986), no 3, 431-454.
\bibitem{STr} Struwe, M., On the evolution of harmonic mappings of Riemannian surfaces, Comment. Math.
Helv. 60 (1985), no. 4, 558-581.

\end{thebibliography}
\end{document}